\documentclass[10 pt]{article}
\usepackage{pgf,tikz}
\usetikzlibrary{arrows}
\usepackage{amsmath,amsfonts,amssymb,amsthm}
\usepackage{graphicx,graphics,epsfig}
\usepackage{epstopdf}	
\usepackage{hyperref}
\usepackage[english]{babel}
\usepackage{makeidx}
\usepackage{latexsym}
\usepackage[latin1]{inputenc}
\newtheorem{thm}{Theorem}[section]

\newtheorem{lem}[thm]{Lemma}

\theoremstyle{definition}
\newtheorem*{rmk}{Remark}
\newtheorem*{ackn}{Acknowledgments}
\DeclareMathOperator{\sgn}{sgn} 
\begin{document} 
%
%
\baselineskip=1.0\baselineskip
\title{\textbf{Sharp Well-Posedness Results for the\\ 
Schr\"odinger-Benjamin-Ono System}}
\author{Leandro Domingues}
\date{December 16, 2014} 
\maketitle \large{\vspace{-0.5cm} 
{\scriptsize \centerline{Departamento de Matemática Aplicada, CEUNES/UFES}
 \centerline{Rodovia BR 101 Norte, Km 60, Bairro Litorâneo, CEP 29932-540, São Mateus, ES, Brazil.}
 \centerline{email: leandro.domingues@ufes.br, leandro.ceunes@gmail.com}} 
\vspace{0.5cm} 
\begin{abstract} 
This work is concerned with the 
Cauchy problem for a coupled Schrödinger-Benjamin-Ono system 
\begin{equation*}
\left \{
    \begin{array}{l}
\vspace{0.1cm}
i\partial_tu+\partial_x^2u=\alpha uv, \hfill t\!\in\![-T,T], \ x\!\in\!\mathbb R,\\
\partial_tv+\nu\mathcal H\partial^2_xv=\beta \partial_x(|u|^2),\\
u(0,x)=\phi, \ v(0,x)=\psi, \ \ \ \ \ \ \ \ \ \ \ \ \ \hfill (\phi,\psi)\!\in\!H^{s}(\mathbb R)\!\times\!H^{s'}\!(\mathbb R).
\end{array}
   \right.
\end{equation*}
In the \emph{non-resonant} case $(|\nu|\ne1)$, we prove local well-posedness 
for a large class of initial data.
This improves the results obtained by Bekiranov, Ogawa and Ponce (1998). 
Moreover, 
we prove 
$C^2$\emph{-ill-posedness} at \emph{low-regularity}, and also when the difference of regularity between 
the initial data is large enough. 
As far as we know, this last ill-posedness result is the first of this kind for a nonlinear dispersive system. 
Finally, we also prove that the local well-posedness result obtained by Pecher (2006) in the \emph{resonant} case 
$(|\nu|=1)$ is sharp except for the end-point.
\end{abstract} 
%
%
%
\begin{center}
\item\section{Introduction}
\end{center} 
In \cite{funa}, Funakoshi and Oikawa deduced the following 
Schrödinger-Benjamin-Ono system, 
\begin{equation}\label{SBOsys}
\left \{
    \begin{array}{l}
\vspace{0.2cm} 
i\partial_tu+\partial_x^2u=\alpha uv, \hfill t\!\in\![-T,T], \ x\!\in\!\mathbb R,\\
\vspace{0.1cm} 
\partial_tv+\nu\mathcal H\partial^2_xv=\beta \partial_x(|u|^2),\\
u(0,x)=\phi, \ v(0,x)=\psi, \ \ \ \ \ \ \ \ \ \ \ \ \ \ \hfill (\phi,\psi)\!\in\!H^{s}(\mathbb R)\!\times\!H^{s'}\!(\mathbb R),
\end{array}
   \right.
\end{equation}
\noindent 
where $\mathcal H$ denotes the Hilbert transform, 
$u=u(t,x)$ is a complex-valued function, $v=v(t,x)$ is a real-valued function and $\alpha,\beta, \nu$ 
are real constants such that $\alpha\beta\ne0$.\\ 
\indent The Schrödinger-Benjamin-Ono system \eqref{SBOsys} describes the motion of two fluids 
with different densities under capillary-gravity waves in a deep water flow. The short surface wave is 
usually described by a Schrödinger type equation and the long internal wave is described by some sort of 
wave equation accompanied by a dispersive term (which is a Benjamin-Ono type equation in this case). 
\\ 
\indent 
The natural function spaces to study the local well-posedness (L.W.P.) of 
this system 
are the Sobolev $H^s\!\!\times\!\!H^{s'}$-type spaces. 
Indeed, for a smooth solution $(u,v)$, the following quantities 
are conserved for every $t\in[-T,T]$ 
\begin{displaymath} 
\left \{ 
\begin{array}{l} 
\vspace{0.3cm} 
\|u(t)\|_2^2, 
\\ 
\vspace{0.2cm} 
Im\int u(t,x)\partial_x\overline{u}(t,x)dx+\frac{\alpha}{2\beta}\|v(t)\|_2^2, 
\\ 
\|\partial_xu(t)\|_2^2+\alpha\int v(t,x)|u(t,x)|^2dx-\frac{\alpha\nu}{2\beta}\|D_x^{1/2}v(t)\|_2^2, 
\end{array}
   \right.
\end{displaymath}
where $D_x=\mathcal H\partial_x$. 
\\ 
\indent 
For $|\nu|\ne1$, the \emph{non-resonant} case, 
Bekiranov, Ogawa and Ponce proved in \cite{bop} the L.W.P. of the system 
\eqref{SBOsys}, 
for $(s,s')$ in the half-line 
\begin{displaymath}
\ell:=\left\{(s,s')\in \mathbb R^2\ : \ s'=s-1/2,\ s\ge0\right\}.
\end{displaymath} 
In \cite{AMP}, Angulo, Matheus and Pilod obtained global well-posedness (G.W.P.), 
also for $(s,s')\in\ell$, by using an idea of Colliander, Holmer and Tzirakis \cite{CHT}.\\ 
%
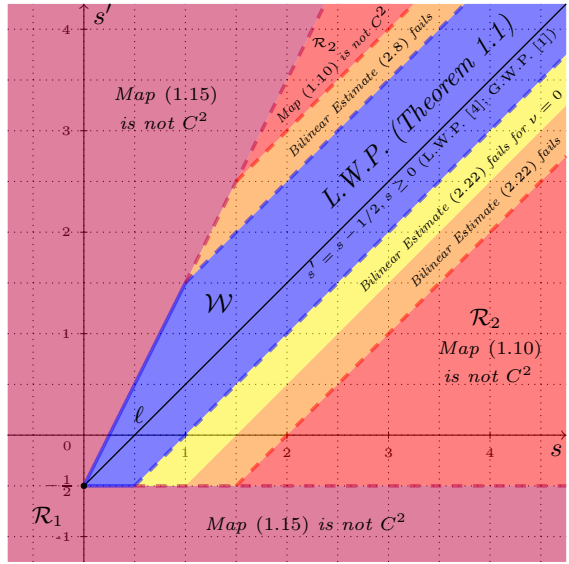
\begin{figure}[!ht] 
\definecolor{cinza}{rgb}{0.7,0.7,0.7}
\begin{center}
\begin{tikzpicture}[scale=1.35][line cap=round,line join=round,>=triangle 45,x=1.0cm,y=1.0cm]
\draw [dotted, xstep=0.5cm,ystep=0.5cm] (-0.75,-1.25) grid (4.75,4.25);
\draw[->] (-0.75,0) -- (4.75,0);
\foreach \x in {1,2,3,4}
\draw[shift={(\x,0)},color=black] (0pt,1pt) -- (0pt,-1pt) node[below] {\tiny\x};
\draw[->] (0,-1.25) -- (0,4.25);
\foreach \y in {-1,1,2,3,4}
\draw[shift={(0,\y)},color=black] (1pt,0pt) -- (-1pt,0pt) node[left] {\tiny\y};
\node at (-0.3,-0.1)[right]{\tiny 0};
\node at (-0.475,-0.5)[right]{\tiny $-\frac12$};
\fill[fill=blue,fill opacity=0.5] (0,-0.5) -- (0.5,-0.5) -- (4.75,3.75) -- (4.75,4.25) -- (3.75,4.25) -- (1,1.5) -- cycle;
\fill[fill=purple,fill opacity=0.5] (-0.75,-1.25) -- (4.75,-1.25) -- (4.75,-0.5) -- 
(0,-0.5) -- (2.375,4.25) -- (-0.75,4.25) -- cycle;
\fill[fill=red,fill opacity=0.5] (1.5,-0.5) -- (4.75,-0.5) -- (4.75,2.75) -- cycle;
\fill[fill=red,fill opacity=0.5] (1.5,2.5) -- (3.25,4.25) -- (2.375,4.25) -- cycle;
\fill[fill=orange,fill opacity=0.5] (1,1.5) -- (3.75,4.25) -- (3.25,4.25) -- (1.5,2.5) -- cycle;
\fill[fill=orange,fill opacity=0.5] (1,-0.5) -- (1.5,-0.5) -- (4.75,2.75) -- (4.75,3.25) -- cycle;
\fill[fill=yellow,fill opacity=0.5] (0.5,-0.5) -- (1,-0.5) -- (4.75,3.25) -- (4.75,3.75) -- cycle;
\draw [line width=0.5pt] (0,-0.5)-- (4.75,4.25);
\draw [color=blue, opacity=0.5][line width=1.35](0,-0.5)-- (0.5,-0.5);
\draw [color=blue, opacity=0.5][line width=1.35](0,-0.5)-- (1,1.5);
\draw [color=blue, opacity=0.5][line width=1.35][dash pattern=on 4pt off 4pt] (0.5,-0.5)-- (4.75,3.75);
\draw [color=blue, opacity=0.5][line width=1.35][dash pattern=on 4pt off 4pt] (1,1.5)-- (3.75,4.25);
\draw [color=purple, opacity=0.5][line width=1.35][dash pattern=on 4pt off 4pt] (0.5,-0.5)-- (4.75,-0.5);
\draw [color=red, opacity=0.5][line width=1.35][dash pattern=on 4pt off 4pt] (1.5,-0.5)-- (4.75,2.75);
\draw [color=purple, opacity=0.5][line width=1.35][dash pattern=on 4pt off 4pt] (1,1.5)-- (2.375,4.25);
\draw [color=red, opacity=0.5][line width=1.35][dash pattern=on 4pt off 4pt] (1.5,2.5)-- (3.25,4.25);
\node at (0.2,3.35)[right]{\scriptsize{\textit{Map \eqref{fluxot}}}};
\node at (0.26,3.1)[right]{\scriptsize{\textit{is not $C^2$}}};
\node at (1.1,-0.875)[right]{\scriptsize{\textit{Map \eqref{fluxot} is not $C^2$}}};
\node at (3.39,0.85)[right]{\scriptsize{\textit{Map \eqref{fluxo}}}};
\node at (3.45,0.6)[right]{\scriptsize{\textit{is not $C^2$}}};
\draw (1.75,3.65) node[right]{$\rotatebox{45}{\tiny{\textit{Map \eqref{fluxo} is not $C^2$}}}$};
\draw (3.1,2.2) node[right]{$\rotatebox{45}{\tiny{\textit{Bilinear Estimate \eqref{NE2} fails}}}$};
\draw (2.6,2.4) node[right]{$\rotatebox{45}{\tiny{\textit{Bilinear Estimate \eqref{NE2} fails for $\nu=0$}}}$};
\draw (1.9,3.45) node[right]{$\rotatebox{45}{\tiny{\textit{Bilinear Estimate \eqref{NE1} fails}}}$};
\node at (1.1,1.3)[right]{$\mathcal W$};
\node at (-0.6,-0.8)[right]{\small{$\mathcal R_1$}};
\node at (3.71,1.15)[right]{\small{$\mathcal R_2$}};
\node at (2.15,3.85)[right]{\tiny{$\mathcal R_2$}};
\node at (0.4,0.2)[right]{\small{$\ell$}};
\draw (2.2,3.175) node[right]{$\rotatebox{45}{\textit{L.W.P. (Theorem \ref{SBO})}}$}; 
\draw (2.05,2.78) node[right]{$\rotatebox{45}{\tiny{$s'=s-1/2, s\ge0$ (L.W.P. \cite{bop};  G.W.P. \cite{AMP})}}$};
\draw [fill] (0,-0.5) circle (0.75pt);
\draw (4.5,-0.3) node [anchor=south west] {$s$};
\draw (0,4.15) node [anchor=west] {$s'$};
\end{tikzpicture}
\caption{Problem \eqref{SBOsys} for $|\nu|\ne1$ and $(\phi,\psi)\!\in\!H^s\!\times\!H^{s'}.$}\label{fig}
\end{center}
\end{figure}
%
\indent 
In Theorem \ref{SBO} of the present paper, we prove the L.W.P. of \eqref{SBOsys} 
for $(s,s')$ in the region 
\begin{equation*}
\mathcal W:=
\{(s,s')\in\mathbb R^2 : -1/2<s'-(s-1/2)<1,\ -1/2 \le s' \le 2s-1/2\}.
\end{equation*} 
Moreover, we establish $C^2$\emph{-ill-posedness} of \eqref{SBOsys} 
for $(s,s')$ 
in the regions 
\begin{equation*}
\mathcal R_1 
:= 
\{(s,s')\in\mathbb R^2 : s'<-1/2\ \ \text{or}\ \ 2s-1/2<s' \} 
\end{equation*} 
and 
\begin{equation*}  
\mathcal R_2 
:= 
\{(s,s')\in\mathbb R^2 : |s'-(s-1/2)|>3/2 \}. 
\end{equation*}
Actually, the ill-posedness result holds in a slightly stronger sense in the region $\mathcal R_1$
(see Theorem \ref{C2ill} for the precise statement). 
Furthermore, Theorem \ref{thm notNE} states that the bilinear estimates used to prove Theorem \ref{SBO} fails 
in a part of the remaining region. 
For the case $\nu=0$, it fails in the entire remaining region. 
All these results are summarized in Figure \ref{fig}. 
In particular, we observe that our results are sharp$^1$ at \emph{low-regularity}. 
\footnotetext[1]{Sharp in the sense that one can not improves the result by performing a Picard iteration, 
since this method provides an analytic flow map data-solution (and hence $C^\infty$ Fréchet differentiable).} 
\\ 
\indent 
For $|\nu|=1$, the \emph{resonant} case, Pecher showed in \cite{per} the L.W.P. of the system 
\eqref{SBOsys} for $(s,s')\in\ell$, except for the end-point $(0,-1/2)$. 
In the present paper, we prove in Theorem \ref{C2ill2} the $C^2$\emph{-ill-posedness} of \eqref{SBOsys} 
for $(s,s')\!\notin\!\ell$. Furthermore, we prove in Theorem \ref{thm endpoint} that the 
key bilinear estimate of Pecher's proof fails at the end-point.\\ 
\indent 
Bekiranov, Ogawa and Ponce also obtained L.W.P. 
for other nonlinear dispersive systems such as 
the Schrödinger-Korteweg-de Vries system (in \cite{bop2}) and the Benney system (in \cite{bop}), 
in both cases, for initial data in $H^s\!\times\!H^{s'}$ with $(s,s')\!\in\!\ell$. 
For the last system, due to scaling properties, the L.W.P. was only investigated in 
$H^{s}(\mathbb R)\!\times\!H^{s-\frac12}\!(\mathbb R)$ (see Remark 2 in \cite{bop}). 
In the case of the system \eqref{SBOsys}, one can scale a solution $(u,v)$ as 
$u_\lambda(t,x)=\lambda^{\frac32}u(\lambda^2t,\lambda x)$, 
$v_\lambda(t,x)=\lambda^2v(\lambda^2t,\lambda x)$. Then $(u_\lambda,v_\lambda)$ solves \eqref{SBOsys} 
with initial data $\phi_\lambda(x)=\lambda^{\frac32}\phi(\lambda x)$ and 
$\psi_\lambda(x)=\lambda^2\psi(\lambda x)$ satisfying 
$\|\phi_\lambda\|_{\dot{H}^s}=\lambda^{1+s}\|\phi\|_{\dot{H}^s}$ and 
$\|\psi_\lambda\|_{\dot{H}^{s'}}=\lambda^{\frac32+s'}\|\psi\|_{\dot{H}^{s'}}$. Thus $s'\!=\!s\!-\!1/2$ 
keeps each norm equivalent under scaling. 
However, Theorem \ref{SBO} shows that the regime $s'=s-1/2$ is not necessary for the L.W.P. 
of the system \eqref{SBOsys}. 
Also, note that Theorem \ref{C2ill} establishes $C^2$\emph{-ill-posedness} 
for $(s,s')$ in a neighborhood of $(-1,-3/2)$, which is a point of critical regularity, 
in the sense that the scaling transformation leaves the 
$\dot{H}^s\!\times\!\dot{H}^{s'}$-norm invariant at this regularity.\\ 
\indent In \cite{GTV}, Ginibre, Tsutsumi and Velo proved the L.W.P. of the Benney system and of the 
1D Zakharov system, for the region 
$\{(s,s')\in\mathbb R^2 : -1/2<s-s'\le1,\ 0 \le s' +1/2 \le 2s\}$. 
In \cite{adan}, Corcho and Linares proved the L.W.P. of the Schrödinger-Korteweg-de Vries system, 
for a region containing the half-line $\ell$. 
Ill-posedness was not investigated in all these works (\cite{bop2}, \cite{GTV}, \cite{bop}, \cite{adan}).\\ 
\indent 
In \cite{wu}, Wu improved the L.W.P. of the Schrödinger-Korteweg-de Vries system obtained in \cite{adan} 
to a larger region. 
Furthermore, he also proved $C^2$\emph{-ill-posedness} results. 
In particular, he showed that his L.W.P. result is sharp$^1$ at \emph{low-regularity}. 
\vspace{0.2cm} \\ 
\indent To state our results, we introduce the integral equations associated to the system \eqref{SBOsys}, 
\begin{eqnarray}
u(t)&=&e^{it\partial^2_x}\phi-i\alpha\int^t_0e^{i(t-t')\partial^2_x}\left(u(t')\cdot v(t')\right)dt',\label{uDuhamel}\\
v(t)&=&e^{-\nu t\mathcal H\partial^2_x}\psi+\beta\int^t_0
e^{-\nu(t-t')\mathcal H\partial^2_x}(\partial_x|u(t')|^2)dt'\label{vDuhamel}, 
\end{eqnarray}
where $e^{it\partial^2_x}$ and $e^{-\nu t\mathcal H\partial^2_x}$ denote the unitary operators for the linear 
Schr\"odinger and Benjamin-Ono equations respectively. We need also to introduce the Bourgain spaces 
for constructing the local solutions. 
For $s,b,s',b',\nu\in\mathbb R$, we let $X^{s,b}$ and $Y^{s',b'}_\nu$ be the completion of the Schwartz class 
$\mathcal S(\mathbb R^2)$ under the norms 
\begin{eqnarray}
\|f\|_{X^{s,b}}
&:=&
\|\langle\xi\rangle^{s}\langle\tau+\xi^2\rangle^{b}\widehat{f}(\tau,\xi)\|_{L^2_{\tau,\xi}}
= 
\|e^{-it\partial^2_x}f\|_{H^b_t(\mathbb R;H^s_x)}\ ,\label{defXsb}\\ 
\|g\|_{Y^{s',b'}_\nu}
&:=&
\|\langle\xi\rangle^{s'}\langle\tau+\nu|\xi|\xi\rangle^{b'}\widehat{g}(\tau,\xi)\|_{L^2_{\tau,\xi}}
= 
\|e^{-\nu t\mathcal H\partial^2_x}g\|_{H^{b'}_t(\mathbb R;H^{s'}_x)}\ ,\label{defYs'b'}  
\end{eqnarray}
where $\langle\cdot\rangle:=\sqrt{1+|\cdot|^2}\ $ and $\ \widehat{f}\ $ is the Fourier transform of $\ f\ $ in 
both $\ x\ $ and $\ t\ $ variables 
\begin{equation*} 
\widehat{f}(\tau,\xi):=
\iint e^{-2\pi i(t\tau+x\xi)}f(x,t)dtdx.
\end{equation*} 
Hereafter, we will simply denote $Y^{s',b'}$ instead of $Y^{s',b'}_\nu$.\\ 
\indent Let $b,b'>1/2$, the Sobolev lemma implies that 
\begin{eqnarray}
X^{s,b} &\hookrightarrow& C^0(\mathbb R;H^s(\mathbb R)),\label{XSobolev} \\ 
Y^{s',b'} &\hookrightarrow& C^0(\mathbb R;H^{s'}(\mathbb R)).\label{YSobolev}
\end{eqnarray} 
Thus, for an interval $I$, 
$M_I:=\{f\in X^{s,b}\ :\ f(t)=0,\ \forall t\in I\}$ is a closed subspace of $X^{s,b}$. 
We define $X^{s,b}_I$ to be the quotient space 
$X^{s,b}/M_I$, which is a Banach space with the norm 
\begin{equation*}
\|f\|_{X^{s,b}_I}:=\inf\{\|\tilde f\|_{X^{s,b}}\ :\ \tilde f(t)=f(t),\ \forall t\in I\}\ . 
\end{equation*}
We write $X^{s,b}_T$ for $X^{s,b}_I$, when $I=[-T,T]$. We define $Y^{s',b'}_T$ similarly. 
\vspace{0.2cm} \\ 
\indent 
Now we are ready to enunciate our results. The first theorem states the L.W.P. of the system 
\eqref{SBOsys}, in the \emph{non-resonant} case, for $(s,s')\in\mathcal W$ (see Figure \ref{fig}). 
\vspace{0.2cm} 
%
\begin{thm}\label{SBO}
Let $|\nu|\ne1$ and $s,s'\in\mathbb R$ satisfying 
\begin{equation}\label{ss'1}
-1/2 \le s' \le 2s-1/2, 
\end{equation}
\begin{equation}\label{ss'2}
s-1 < s' <  s+1/2.
\end{equation}
The Cauchy problem \eqref{SBOsys} is locally well-posed in 
$H^s(\mathbb R)\times H^{s'}(\mathbb R)$, in the following sense:\\ 
For every $R>0$, there exist $T\!=\!T(R)>0$ and $b,b'>1/2$ such that if $\|\phi\|_{H^s}\!+\!\|\psi\|_{H^{s'}}\!<\!R$, 
there exists a unique solution $(u,v)\!\in\! X^{s,b}_T\!\times\! Y^{s',b'}_T$ satisfying 
\eqref{uDuhamel}-\eqref{vDuhamel} for all $t\!\in\![-T,T]$. 
Moreover, this solution satisfies 
\begin{displaymath} 
(u,v)\in C^0([-T,T];H^s(\mathbb R))\times C^0([-T,T];H^{s'}(\mathbb R)), 
\end{displaymath} 
and the associated flow map data-solution,
\begin{equation}\label{fluxo}
S:B_R\to  C^0([-T,T];H^s(\mathbb R))\times C^0([-T,T];H^{s'}(\mathbb R)), \ \ 
(\phi,\psi)\mapsto (u,v),
\end{equation}
is Lipschitz continuous, where $B_R$ is the open ball in $H^s(\mathbb R)\times H^{s'}(\mathbb R)$, centered at 
the origin with radius $R$.
\end{thm}
%
\indent 
Next, we give the main ingredients in the proof of Theorem \ref{SBOsys}. 
Following the procedure employed in \cite{bop}, we use the Banach Fixed Point theorem and the Fourier 
restriction norm method introduced by Bourgain in \cite{Bo}. 
So the difficulty is to extend the following bilinear estimates found in \cite{bop} 
\begin{equation}\label{NE1bop} 
\|\partial_x(u_1\overline{u_2})\|_{Y^{s-\frac12,a}}
\le 
C\|u_1\|_{X^{s,b}}\|u_2\|_{X^{s,b}}\ ,
\hspace{1.35cm} 
b>1/2,\ a\le0,\ s\ge0, 
\end{equation} 
\begin{equation}\label{NE2bop} 
\|uv\|_{X^{s,a}} 
\le 
C\|u\|_{X^{s,b}}\|v\|_{Y^{s-\frac12,b}}\ , 
\hspace{0.9cm} 
3/4>b>1/2,\ a<-1/4,\ s\ge0, 
\end{equation} 
to new ones. 
Proceeding as in \cite{GTV}, 
we decouple the modulation regularities of the spaces $X$ and $Y$ in order to gain spatial regularity 
(i.e., we replace $(s\!-\!1/2,b)$ by $(s',b')$ in $Y$). 
Then, by choosing $1/2\!<\!b\!<\!c\!<\!3/4$ and $1/2\!<\!b'\!<\!c'\!<\!3/4$ depending on $(s,s')$, 
we prove the following estimates  
\begin{equation}\label{NE1intro} 
\|\partial_x(u_1\overline{u_2})\|_{Y^{s',c'-1}}
\le
C\|u_1\|_{X^{s,b}}\|u_2\|_{X^{s,b}} \ , 
\end{equation} 
%
\begin{equation}\label{NE2intro} 
\|uv\|_{X^{s,c-1}} 
\le 
C\|u\|_{X^{s,b}}\|v\|_{Y^{s',b'}} \ , 
\end{equation} 
for $(s,s')$ in larger regions (c.f. Theorems \ref{thm NE1} and \ref{thm NE2}). 
Hence, the system \eqref{SBOsys} is L.W.P. for $(s,s')\in\mathcal W$, 
where both estimates \eqref{NE1intro} and \eqref{NE2intro} hold. 
The estimate \eqref{NE1bop} offers minor difficulty in \cite{bop}, 
since the regime $s'\!=\!s\!-\!1/2,\ s\ge0$ 
allows good cancellations in the frequency interactions. 
However, those cancellations do not occur anymore 
for $(s,s')$ in the larger region where the estimate \eqref{NE1intro} holds. 
Thus, we need to perform a new decomposition of the Euclidean space 
(c.f. \eqref{regANE1}-\eqref{regB2NE1}) in order to obtain \eqref{NE1intro}. 
On the other hand, there are no good cancellations for the estimate \eqref{NE2bop}, 
even in the regime $s'\!=\!s\!-\!1/2,\ s\ge0$. 
However, 
we are able to prove the estimate \eqref{NE2intro} for $(s,s')$ in a larger region 
by performing the decomposition \eqref{regANE2}-\eqref{regC2NE2}, 
which is slightly different from the one used in \cite{bop} to obtain the estimate \eqref{NE2bop}. 
\vspace{0.2cm} 
\\ 
\indent 
In the next theorem, we state an ill-posedness result for the \emph{non-resonant} case. 
%
\begin{thm}\label{C2ill}
Let $|\nu|\ne1$ and $s,s'\in\mathbb R$. 
Suppose that the Cauchy problem \eqref{SBOsys} is locally well-posed in 
$H^{s}(\mathbb R)\times H^{s'}(\mathbb R)$, in the sense of Theorem \ref{SBO}. 
\begin{itemize}
\item[(i)] 
If \eqref{ss'1} is not verified, the associated flow map data-solution, 
\begin{equation}\label{fluxot}
S^t: B_R\to H^s(\mathbb R)\times H^{s'}(\mathbb R), \ \ (\phi,\psi)\mapsto (u(t),v(t)),
\end{equation}
is not $C^2$ at zero$^1$, for $t\in [-T,0)\cup(0,T]$. 
Neither is, a fortiori, the flow map \eqref{fluxo}.
\footnotetext[1]
{Actually, we prove that these maps are not two times Fréchet differentiable at zero.} 
\item[(ii)] 
If $\left|s'-(s-1/2)\right|>3/2$, the associated map data-solution \eqref{fluxo} is not $C^2$ at zero$^1$. 
\end{itemize}
\vspace{0.2cm} 
\end{thm} 
%
The first \emph{$C^2$-ill-posedness} result of this kind was proved by Tzvetkov in 
\cite{Tzvet} for the KdV equation. 
We essentially follow his argument to prove Theorem \ref{C2ill} $(i)$. 
There is an additional technical difficulty to prove $(ii)$. 
To overcome this difficulty, we allow the variable $t$ to move. 
Therefore, $(ii)$ presents a conclusion for the flow map \eqref{fluxo} instead of the flow map \eqref{fluxot}. 
We emphasize that this approach has already been used in previous works (e.g., \cite{BTao} and \cite{LG}).
\begin{rmk} 
As far as we know, Theorem \ref{C2ill} $(ii)$ is the first result concerning the ill-posedness of a nonlinear 
dispersive system when the difference of regularity between the initial data is large enough 
(see region $\mathcal R_2$ in Figure \ref{fig}). 
Such result seems natural, due to the coupling of the system via the nonlinearities. 
We believe that the same approach used to prove Theorem \ref{C2ill} $(ii)$ can provide similar results for other 
nonlinear dispersive systems such as the Zakharov system and the Schrödinger-Korteweg-de Vries system. 
We plan to address this issue in a forthcoming paper. 
\end{rmk} 
Finally, we state an ill-posedness result for the \emph{resonant} case. 
%
\begin{thm}\label{C2ill2}
Let $|\nu|=1$ and $(s,s')\notin\ell$, i.e., $s'\ne s-1/2$ or $s<0$. 
If the Cauchy problem \eqref{SBOsys} is locally 
well-posed in $H^{s}(\mathbb R)\times H^{s'}(\mathbb R)$, 
the flow map data-solution \eqref{fluxot} is not $C^2$ at zero$^1$, for $t\in [-T,0)\cup(0,T]$ 
and, a fortiori, 
neither is the flow map \eqref{fluxo}.
\end{thm} 
\vspace{0.2cm} 
\indent 
Throughout the whole text, we use the following notations:
\begin{itemize} 
\item For any $x\in\mathbb R$, we define 
$\sgn(x):=x{|x|}^{-1}$ if $x\ne0$ and $\sgn(x):=1$ if $x=0$. 
\item Let $\mathbf 1_\Omega$ denotes the characteristic function of an arbitrary set $\Omega$, 
i.e., 
$\mathbf 1_\Omega(x)=1$ if $x\in\Omega$ and $\mathbf 1_\Omega(x)=0$ if $x\notin\Omega$. 
\item Fix $\eta$ a smooth function supported on the
interval $[-2,2]$ such that $\eta(x)\equiv 1$ for all $|x|\leq 1$
and, for each $T>0$, $\eta_T(t):=\eta(t/T)$. 
\item For positive quantities $X$ and $Y$, the notation $X \lesssim Y$ means that there exist a constant 
$C>0$ such that $X\le CY$, depending only on 
the parameters $\alpha, \beta$ and $\nu$ related to (\ref{SBOsys}), 
on the indices $s, s', b,c,b'$ and $c'$ related to the Bourgain spaces 
in the bilinear estimates (\ref{NE1}) and (\ref{NE2}), and  
on certain norms of the fixed cut-off function $\eta$. 
We denote $X\gtrsim Y$ when $Y\lesssim X$, and denote $X\sim Y$ when $X \lesssim Y\lesssim X$. 
\end{itemize} 
\ 
\indent 
This paper is organized as follows. In Section 2, we establish the new bilinear estimates that we use to prove 
Theorem \ref{SBO} in Section 3. In Section 4, we prove Theorems \ref{C2ill} , \ref{C2ill2}, \ref{thm notNE} and \ref{thm endpoint}. 
%
%
%
\begin{center}
\item\section{Bilinear Estimates}
\end{center}
\setcounter{section}{2} \setcounter{equation}{0} 
\indent 
In this section, we improve the bilinear estimates presented in \cite{bop}. First, we state some 
calculus inequalities which will be useful in the proofs of Theorems \ref{thm NE1} and \ref{thm NE2}.\\ 
%
\begin{lem} \label{lem CalcElem}
Let $\alpha>1/2$ and $1/2<\beta,\gamma\le1$. Then, for all $p\ne0$ and $q,r\in \mathbb R$, 
the following estimates hold:
\begin{eqnarray}
(i)&&\int\frac{dx}{{\langle x-q\rangle}^{2\beta}{\langle x-r\rangle}^{2\gamma}}
\lesssim\frac{1}{{\langle q-r\rangle}^{2\min\{\beta,\gamma\}}}\ ,\label{CEb}\\
\nonumber\\
(ii)&&\int\frac{dx}{{\langle x-q\rangle}^{2\beta}{\langle x-r\rangle}^{2(1-\gamma)}}
\lesssim \frac{1}{{\langle q-r\rangle}^{2(1-\gamma)}},\label{CE1-b}\ \\
\nonumber\\
(iii)&&\int\frac{dx}{{\langle px^2+qx+r\rangle}^\alpha}\lesssim \frac{1}{|p|}.\hspace{10.0cm}\label{CEQ}
\end{eqnarray}
%
\end{lem}
\noindent 
The estimates \eqref{CEb} and \eqref{CE1-b} are particular cases of the estimates established in 
Lemma 4.2 of \cite{GTV}.
The estimate \eqref{CEQ} follows from elementary computations (for the ideas, see (2.14) of \cite{bop2} and 
note that $\langle \cdot\rangle\sim1+|\cdot|$). 
%
\begin{thm} \label{thm NE1}
Assume that $|\nu|\ne1$. Let $s,s'\in\mathbb R$ be such that $s\ge0$,
\begin{eqnarray}
s' &\le& 2s-1/2,\label{ss'1r}\\
s'&<&s+1/2.\label{ss'2r}
\end{eqnarray}
Then, for all $b,c'\in \mathbb R$ such that
\begin{equation}\label{NE1b}
\max\left\{1/2\ , \ (s'-s)/2+1/2\right\}<b,
\end{equation}
\begin{equation}\label{NE1c'}
c'<\min\left\{3/4-(s'-s)/2\ , \ 3/4\right\},
\end{equation}
the following estimate holds:
\begin{equation}\label{NE1}
\|\partial_x(u_1\overline{u_2})\|_{Y^{s',c'-1}}
\lesssim
\|u_1\|_{X^{s,b}}\|u_2\|_{X^{s,b}} \ , \ \ \ \forall u_1,u_2\in X^{s,b}.
\end{equation}
\end{thm} 
%
%
\noindent \textbf{Proof.} 
It is sufficient to show \eqref{NE1} for $u_1,u_2\in \mathcal S(\mathbb R^2)$. Thus, letting 
\begin{displaymath}
f(\tau,\xi):={\langle\xi\rangle}^s{\langle\tau+\xi^2\rangle}^b\ \widehat{u_1}(\tau,\xi), \ \
g(\tau,\xi):={\langle\xi\rangle}^s{\langle\tau-\xi^2\rangle}^b\ \overline{\widehat{u_2}}(-\tau,-\xi), 
\end{displaymath} 
and denoting $\tau_2:=\tau -\tau_1$ and $\xi_2:=\xi -\xi_1$, the estimate (\ref{NE1}) is equivalent to
\begin{displaymath}
\left\|\frac{i\xi{\langle\xi\rangle}^{s'}}{{\langle \tau+\nu|\xi|\xi\rangle}^{1-c'}}
\iint \frac{f(\tau_2,\xi_2)g(\tau_1,\xi_1)d\tau_1d\xi_1}
{{\langle\xi_2\rangle}^s{\langle\tau_2+\xi_2^2\rangle}^b{\langle\xi_1\rangle}^s{\langle\tau_1-\xi_1^2\rangle}^b}
\right\|_{L^2_{\tau,\xi}}
\lesssim
\|f\|_{L^2}\|g\|_{L^2},\ \forall f,g\in\mathcal S(\mathbb R^2).
\end{displaymath}
For convenience, we rewrite this estimate as
\begin{equation}\label{NE1eqv}
\left\|\iint\Phi(\tau,\xi,\tau_1,\xi_1)f(\tau_2,\xi_2)g(\tau_1,\xi_1)d\tau_1d\xi_1\right\|_{L^2_{\tau,\xi}}
\lesssim \|f\|_{L^2}\|g\|_{L^2},\ \ 
\forall f,g\in\mathcal S(\mathbb R^2),
\end{equation}
where
$$\Phi(\tau,\xi,\tau_1,\xi_1):=\frac{i\xi{\langle\xi\rangle}^{s'}{\langle\sigma\rangle}^{c'-1}}
{{\langle\xi_2\rangle}^s{\langle\sigma_2\rangle}^b{\langle\xi_1\rangle}^s{\langle\sigma_1\rangle}^b} \ ,$$
with the additional notation
$\sigma:=\tau+\nu|\xi|\xi$, $\sigma_1:=\tau_1-\xi_1^2$ 
and $\sigma_2:=\tau_2+\xi_2^2.$ With this notation, 
the algebraic relation associated to (\ref{NE1eqv}) is given by
\begin{equation}\label{algrel}
\sigma-\sigma_1-\sigma_2=2\xi\xi_1-(1-\nu\sgn(\xi))\xi^2=(1+\nu\sgn(\xi))\xi^2-2\xi\xi_2.
\end{equation}
We split $\mathbb R^4$ into the following regions
\begin{eqnarray}
\mathcal{A}_{ \ } &=& \left\{(\tau,\xi,\tau_1,\xi_1) \in \mathbb{R}^4 \ : \
|(1-\nu\sgn(\xi))\xi-2\xi_1| < c_\nu|\xi| \right\},\label{regANE1}   \\
\mathcal{B}_{ \ } &=& \left\{(\tau,\xi,\tau_1,\xi_1) \in \mathcal{A}^c\ :
|\sigma|=\max\{|\sigma|,|\sigma_1|,|\sigma_2|\}\right\},\label{regBNE1} \\
\mathcal{B}_{1} &=& \left\{(\tau,\xi,\tau_1,\xi_1) \in \mathcal{A}^c\ :
|\sigma_1|=\max\{|\sigma|,|\sigma_1|,|\sigma_2|\}\right\},\label{regB1NE1} \\
\mathcal{B}_{2} &=& \left\{(\tau,\xi,\tau_1,\xi_1) \in \mathcal{A}^c\ :
|\sigma_2|=\max\{|\sigma|,|\sigma_1|,|\sigma_2|\}\right\},\label{regB2NE1} 
\end{eqnarray}
where $c_\nu:=\frac{|1-|\nu||}{2}>0.$\\
We use the Cauchy-Schwarz inequality and the Fubini theorem to estimate the left-hand side of (\ref{NE1eqv}) 
restricted to each one of these sets (also perform a change of variables in the region $\mathcal{B}_{2}$). 
Thus (\ref{NE1eqv}) is a consequence of the following estimates 
\begin{eqnarray}
\|\mathbf 1_{\mathcal{A}\cup\mathcal{B}}\Phi\|_{L^\infty_{\tau,\xi}(L^2_{\tau_1,\xi_1})}
&\lesssim& 1,\label{NE1AB}\\
\|\mathbf 1_{\mathcal{B}_1}\Phi\|_{L^\infty_{\tau_1,\xi_1}(L^2_{\tau,\xi})}&\lesssim& 1,\label{NE1B1}\\
\|\mathbf 1_{\widetilde{\mathcal{B}}_2}\widetilde{\Phi}\|_{L^\infty_{\tau_2,\xi_2}(L^2_{\tau,\xi})}
&\lesssim& 1,\label{NE1B2}
\end{eqnarray}
where 
\begin{eqnarray*}
\widetilde{\Phi}(\tau,\xi,\tau_2,\xi_2)&:=&\Phi(\tau,\xi,\tau-\tau_2,\xi-\xi_2),\ 
\forall(\tau,\xi,\tau_2,\xi_2)\in \mathbb R^4,\\
\widetilde{\mathcal{B}}_2&:=&\left\{(\tau,\xi,\tau_2,\xi_2) \in \mathbb R^4\ :
(\tau,\xi,\tau-\tau_2,\xi-\xi_2)\in\mathcal{B}_2\right\}.
\end{eqnarray*}
\noindent\textit{Proof of the estimate (\ref{NE1AB})}: In the region $\mathcal{A}$, $|\xi|\sim|\xi_1|\sim|\xi_2|$. 
In fact, rewriting  
\begin{equation*}
2|\xi_1|=|(1-\nu\sgn(\xi))\xi-2\xi_1-(1-\nu\sgn(\xi))\xi|,
\end{equation*}
we conclude that $c_\nu|\xi|\le 2|\xi_1|\le (c_\nu+1+|\nu|)|\xi|$.
Similarly, we have 
\begin{equation*}
2|\xi_2|=|(1-\nu\sgn(\xi))\xi-2\xi_1+(1+\nu\sgn(\xi))\xi|, 
\end{equation*}
thus $c_\nu|\xi|\le 2|\xi_2|\le (c_\nu+1+|\nu|)|\xi|$. 
Hence, we get from $c'<3/4$ and \eqref{ss'1r} that 
\begin{displaymath}
|\Phi(\tau,\xi,\tau_1,\xi_1)|\lesssim 
\frac{{\langle\xi\rangle}^{s'-2s+\frac12}|\xi|^\frac12}{{\langle\sigma_2\rangle}^b{\langle\sigma_1\rangle}^b}
\lesssim 
\frac{|\xi|^\frac12}{{\langle\sigma_2\rangle}^b{\langle\sigma_1\rangle}^b}, \ 
\forall(\tau,\xi,\tau_1,\xi_1)\in \mathcal{A}.
\end{displaymath}
The same estimate holds in $\mathcal{B}$. In fact, in the region $\mathcal{A}^c$, 
we have from \eqref{algrel} that 
\begin{equation}\label{ximod}
c_{\nu}{|\xi|}^2\le |(1-\nu\sgn(\xi))\xi^2-2\xi\xi_1|=|\sigma-\sigma_1-\sigma_2|.
\end{equation}
In particular, ${|\xi|}^2\lesssim|\sigma|$ in the region $\mathcal{B}$. 
Note also that $\langle\xi_1+\xi_2\rangle\lesssim\langle\xi_1\rangle\langle\xi_2\rangle$. Thus, 
we deduce from \eqref{NE1c'} that 
\begin{displaymath}
|\Phi(\tau,\xi,\tau_1,\xi_1)|\lesssim 
\frac{{\langle\xi\rangle}^{s'-s+\frac12}}{{\langle\sigma\rangle}^{1-c'}}\cdot
\frac{|\xi|^\frac12}{{\langle\sigma_2\rangle}^b{\langle\sigma_1\rangle}^b}\lesssim 
\frac{|\xi|^\frac12}{{\langle\sigma_2\rangle}^b{\langle\sigma_1\rangle}^b}, \ 
\forall(\tau,\xi,\tau_1,\xi_1)\in \mathcal{B}. 
\end{displaymath}
Now, observe from \eqref{CEb} and $b>1/2$, that  
\begin{displaymath}
\sup_{\tau,\xi}\left\|\frac{|\xi|^\frac12}{{\langle\sigma_2\rangle}^b
{\langle\sigma_1\rangle}^b}\right\|_{L^2_{\tau_1,\xi_1}}
\lesssim
\sup_{\tau,\xi}{\left[\int\frac{|\xi|}
{{\langle{2\xi\xi_1-\tau-\xi^2}\rangle}^{2b}}d\xi_1\right]}^\frac12
\lesssim1. 
\end{displaymath}
This concludes the proof of \eqref{NE1AB}.\\ 
%
\textit{Proof of the estimate (\ref{NE1B1})}: By (\ref{ximod}), ${|\xi|}^2\lesssim|\sigma_1|$ in the region 
$\mathcal{B}_1$. Thus \eqref{NE1b} implies that 
\begin{displaymath}
|\Phi(\tau,\xi,\tau_1,\xi_1)|\lesssim 
\frac{{\langle\xi\rangle}^{s'-s+1}}{{\langle\sigma_1\rangle}^{b}}\cdot
\frac{1}{{\langle\sigma_2\rangle}^b{\langle\sigma\rangle}^{1-c'}}\lesssim 
\frac{1}{{\langle\sigma_2\rangle}^b{\langle\sigma\rangle}^{1-c'}}, \ 
\forall(\tau,\xi,\tau_1,\xi_1)\in \mathcal{B}_1. 
\end{displaymath}
From \eqref{CE1-b}, \eqref{algrel}, \eqref{CEQ} and $c'<3/4$, we deduce that 
\begin{displaymath}
\sup_{\tau_1,\xi_1}\left\|\frac{1}{{\langle\sigma_2\rangle}^b
{\langle\sigma\rangle}^{1-c'}}\right\|_{L^2_{\tau,\xi}}
\lesssim
\sup_{\tau_1,\xi_1}{\left[\int\frac{d\xi}
{{\langle{(1-\nu\sgn(\xi))\xi^2-2\xi\xi_1-\sigma_1}\rangle}^{2(1-c')}}\right]}^\frac12
\lesssim1,
\end{displaymath}
which yields \eqref{NE1B1}.\\ 
%
\textit{Proof of the estimate (\ref{NE1B2})}: Denoting $\tau_1:=\tau-\tau_2$, $\xi_1:=\xi-\xi_2$ and 
$\sigma, \ \sigma_1, \ \sigma_2$ as before,  we have ${|\xi|}^2\lesssim|\sigma_2|$ 
in the region $\widetilde{\mathcal{B}}_2$. Then, we deduce from \eqref{NE1b} that 
\begin{displaymath}
|\widetilde{\Phi}(\tau,\xi,\tau_2,\xi_2)|\lesssim 
\frac{{\langle\xi\rangle}^{s'-s+1}}{{\langle\sigma_2\rangle}^{b}}\cdot
\frac{1}{{\langle\sigma_1\rangle}^b{\langle\sigma\rangle}^{1-c'}}\lesssim 
\frac{1}{{\langle\sigma_1\rangle}^b{\langle\sigma\rangle}^{1-c'}}, \ 
\forall(\tau,\xi,\tau_2,\xi_2)\in \widetilde{\mathcal{B}}_2, 
\end{displaymath}
and from \eqref{CE1-b}, \eqref{algrel}, \eqref{CEQ} and $c'<3/4$ that 
\begin{displaymath}
\sup_{\tau_2,\xi_2}\left\|\frac{1}{{\langle\sigma_1\rangle}^b
{\langle\sigma\rangle}^{1-c'}}\right\|_{L^2_{\tau,\xi}}
\lesssim
\sup_{\tau_2,\xi_2}{\left[\int\frac{1}
{{\langle{(1+\nu\sgn(\xi))\xi^2-2\xi\xi_2+\sigma_2}\rangle}^{2(1-c')}}d\xi\right]}^\frac12
\lesssim1,
\end{displaymath}
which concludes \eqref{NE1B2}. This finishes the proof of \eqref{NE1}.\hfill$\square$\\ 
%
\begin{thm} \label{thm NE2}
Assume that $|\nu|\ne1$. Let $s,s'\in\mathbb R$ be such that $s\ge0$,
\begin{eqnarray}
-1/2 &\le& s',\label{ss'1l}\\
s-1&<&s'.\label{ss'2l}
%
%
\end{eqnarray}
Then, for all $b, b', c\in\mathbb R$ such that $1/2<b,b'$ and 
\begin{equation}\label{NE2c}
1/2<c<\min\left\{3/4\ , \ (s'-s)/2+1\right\},
%
%
\end{equation}
the following estimate holds:
\begin{equation}\label{NE2}
\|uv\|_{X^{s,c-1}}
\lesssim
\|u\|_{X^{s,b}}\|v\|_{Y^{s',b'}} \ , \ \ \ \forall u\in X^{s,b}, \ \forall v\in Y^{s',b'}.
\end{equation}
\end{thm}
%
\noindent \textbf{Proof.} 
%
It is sufficient to show \eqref{NE1} for $u,v\in \mathcal S(\mathbb R^2)$. Thus letting 
$$f(\tau,\xi):={\langle\xi\rangle}^s{\langle\tau+\xi^2\rangle}^b\hat u(\tau,\xi), \ \ \
g(\tau,\xi):={\langle\xi\rangle}^{s'}{\langle\tau+\nu|\xi|\xi\rangle}^{b'} \hat{v}(\tau,\xi),$$ 
and denoting $\tau_2:=\tau -\tau_1$ and $\xi_2:=\xi -\xi_1$, the estimate (\ref{NE2}) is equivalent to
\begin{displaymath}
\left\|\frac{{\langle\xi\rangle}^{s}}{{\langle \tau+\xi^2\rangle}^{1-c}}
\iint \frac{f(\tau_2,\xi_2)g(\tau_1,\xi_1)d\tau_1d\xi_1}
{{\langle\xi_2\rangle}^s{\langle\tau_2+\xi_2^2\rangle}^b
{\langle\xi_1\rangle}^{s'}{\langle\tau_1+\nu|\xi_1|\xi_1\rangle}^{b'}}
\right\|_{L^2_{\tau,\xi}}
\lesssim
\|f\|_{L^2}\|g\|_{L^2},\ \forall f,g\in\mathcal S(\mathbb R^2).
\end{displaymath}
For convenience, we rewrite this estimate as
\begin{equation}\label{NE2eqv}
\left\|\iint\Psi(\tau,\xi,\tau_1,\xi_1)f(\tau_2,\xi_2)g(\tau_1,\xi_1)d\tau_1d\xi_1\right\|_{L^2_{\tau,\xi}}
\lesssim \|f\|_{L^2}\|g\|_{L^2},\ \ 
\forall f,g\in\mathcal S(\mathbb R^2).
\end{equation}
where
\begin{displaymath} 
\Psi(\tau,\xi,\tau_1,\xi_1):=\frac{{\langle\xi\rangle}^s{\langle\sigma\rangle}^{c-1}}
{{\langle\xi_2\rangle}^s{\langle\sigma_2\rangle}^b{\langle\xi_1\rangle}^{s'}{\langle\sigma_1\rangle}^{b'}} \ , 
\end{displaymath}
with the additional notation
$\sigma:=\tau+\xi^2$, $\sigma_1:=\tau_1+\nu|\xi_1|\xi_1$ and $\sigma_2:=\tau_2+\xi_2^2$. 
With this notation, the algebraic relation associated to (\ref{NE2eqv}) is given by
\begin{equation}\label{algrel2}
\sigma-\sigma_1-\sigma_2=2\xi\xi_1-(1+\nu\sgn(\xi_1))\xi_1^2=(1-\nu\sgn(\xi_1))\xi_1^2+2\xi_1\xi_2 \ .
\end{equation}
We split $\mathbb R^4$ into the following regions
\begin{eqnarray}
\mathcal{A}_{ \ } &=& \left\{(\tau,\xi,\tau_1,\xi_1) \in \mathbb{R}^4 \ : \ |\xi_1|\le1 \right\},\label{regANE2} \\
\mathcal{B}_{ \ } &=& \left\{(\tau,\xi,\tau_1,\xi_1) \in \mathcal{A}^c\ :
|(1+\nu\sgn(\xi_1))\xi_1-2\xi| < c_\nu|\xi_1|\right\},\label{regBNE2} \\
\mathcal{C}_{ \ } &=& \left\{(\tau,\xi,\tau_1,\xi_1) \in \mathcal{A}^c\cap\mathcal{B}^c\ :
|\sigma|=\max\{|\sigma|,|\sigma_1|,|\sigma_2|\}\right\},\label{regCNE2} \\
\mathcal{C}_{1} &=& \left\{(\tau,\xi,\tau_1,\xi_1) \in \mathcal{A}^c\cap\mathcal{B}^c\ :
|\sigma_1|=\max\{|\sigma|,|\sigma_1|,|\sigma_2|\}\right\},\label{regC1NE2} \\
\mathcal{C}_{2} &=& \left\{(\tau,\xi,\tau_1,\xi_1) \in \mathcal{A}^c\cap\mathcal{B}^c\ :
|\sigma_2|=\max\{|\sigma|,|\sigma_1|,|\sigma_2|\}\right\},\label{regC2NE2} 
\end{eqnarray}
where $c_\nu:=\frac{|1-|\nu||}{2}>0.$\\
Arguing similarly 
to the proof of Theorem \ref{thm NE1}, 
it is enough to show 
\begin{eqnarray}
\|\mathbf 1_{\mathcal{A}\cup\mathcal{B}\cup\mathcal{C}}\Psi\|_{L^\infty_{\tau,\xi}(L^2_{\tau_1,\xi_1})}
&\lesssim& 1,\label{NE2ABC}\\
\|\mathbf 1_{\mathcal{C}_1}\Psi\|_{L^\infty_{\tau_1,\xi_1}(L^2_{\tau,\xi})}&\lesssim& 1,\label{NE2C1}\\
\|\mathbf 1_{\widetilde{\mathcal{C}}_2}\widetilde{\Psi}\|_{L^\infty_{\tau_2,\xi_2}(L^2_{\tau_1,\xi_1})}
&\lesssim& 1,\label{NE2C2}
\end{eqnarray}
where
\begin{eqnarray*}
\widetilde{\Psi}(\tau_2,\xi_2,\tau_1,\xi_1)&:=&\Psi(\tau_1+\tau_2,\xi_1+\xi_2,\tau_1,\xi_1),\ 
\forall(\tau_2,\xi_2,\tau_1,\xi_1)\in \mathbb R^4,\\
\widetilde{\mathcal{C}}_2&:=&\left\{(\tau_2,\xi_2,\tau_1,\xi_1) \in \mathbb R^4\ :
(\tau_1+\tau_2,\xi_1+\xi_2,\tau_1,\xi_1)\in\mathcal{C}_2\right\}.
\end{eqnarray*}
\noindent\textit{Proof of the estimate (\ref{NE2ABC})}: 
In the region $\mathcal{A}$, we get that  
\begin{displaymath}
|\Psi(\tau,\xi,\tau_1,\xi_1)|\lesssim 
\frac{{\langle\xi_1\rangle}^{s-s'}}{{\langle\sigma_2\rangle}^b{\langle\sigma_1\rangle}^{b'}}\lesssim 
\frac{1}{{\langle\sigma_2\rangle}^b{\langle\sigma_1\rangle}^{b'}}, \ 
\forall(\tau,\xi,\tau_1,\xi_1)\in \mathcal{A},
\end{displaymath}
since $c<3/4$ and $|\xi_1|\le1$. Therefore, we deduce from \eqref{CEb}, \eqref{algrel2} 
and \eqref{CEQ} that 
\begin{equation}\label{psiA}
\|\mathbf 1_{\mathcal{A}}\Psi\|_{L^\infty_{\tau,\xi}(L^2_{\tau_1,\xi_1})}
\lesssim
\sup_{\tau,\xi}{\left[\int\frac{d\xi_1}
{{\langle{(1+\nu\sgn(\xi_1))\xi_1^2-2\xi\xi_1+\sigma}\rangle}^{2\min\{b,b'\}}}\right]}^\frac12
\lesssim1.
\end{equation}
In the region $\mathcal{B}$, $|\xi|\lesssim|\xi_2|$. Indeed, the identities 
\begin{eqnarray*}
2|\xi|&=&|(1+\nu\sgn(\xi_1))\xi_1-2\xi-(1+\nu\sgn(\xi_1))\xi_1|\\
2|\xi_2|&=&|(1+\nu\sgn(\xi_1))\xi_1-2\xi+(1-\nu\sgn(\xi_1))\xi_1|
\end{eqnarray*}
imply that $2|\xi|\le (c_\nu+1+|\nu|)|\xi_1|$ and $c_\nu|\xi_1|\le 2|\xi_2|$. 
Therefore,
\begin{equation*}
|\Psi(\tau,\xi,\tau_1,\xi_1)|\lesssim 
\frac{1}{{\langle\xi_1\rangle}^{s'}{\langle\sigma_2\rangle}^b{\langle\sigma_1\rangle}^{b'}}
\lesssim
\frac{{|2(1+\nu\sgn(\xi_1))\xi_1-2\xi|}^\frac12}
{{\langle\xi_1\rangle}^{s'+\frac12}{\langle\sigma_2\rangle}^b{\langle\sigma_1\rangle}^{b'}},\ \ 
\forall(\tau,\xi,\tau_1,\xi_1)\in \mathcal{B},
\end{equation*}
since $\langle\xi_1\rangle\le \sqrt2|\xi_1|\le \frac{\sqrt2}{c_{\nu}}|2(1+\nu\sgn(\xi_1))\xi_1-2\xi|$.
Hence, performing the change of variable $\eta:=(1+\nu\sgn(\xi_1))\xi_1^2-2\xi\xi_1+\sigma$, 
we conclude from \eqref{ss'1l}, \eqref{CEb}, \eqref{algrel2} and $b',b>1/2$,
\begin{equation}\label{psiB}
\|\mathbf 1_{\mathcal{B}}\Psi\|_{L^\infty_{\tau,\xi}(L^2_{\tau_1,\xi_1})}
\lesssim
\sup_{\tau,\xi}{\left[\int\frac{d\eta}{{\langle\eta\rangle}^{2\min\{b,b'\}}}\right]}^\frac12
\lesssim1.
\end{equation}
In the region ${\mathcal{A}}^c\cap{\mathcal{B}}^c$, we have 
\begin{equation}\label{xi1mod}
c_{\nu}{|\xi_1|}^2\le |(1-\nu\sgn(\xi_1))\xi_1^2-2\xi\xi_1|=|\sigma-\sigma_1-\sigma_2|.
\end{equation}
In particular, ${|\xi_1|}^2\lesssim|\sigma|$ in the region $\mathcal{C}$. Thus, \eqref{NE2c} implies 
\begin{displaymath}
|\Psi(\tau,\xi,\tau_1,\xi_1)|\lesssim 
\frac{{\langle\xi_1\rangle}^{s-s'}}{{\langle\sigma\rangle}^{1-c}}\cdot
\frac{1}{{\langle\sigma_2\rangle}^b{\langle\sigma_1\rangle}^{b'}}\lesssim 
\frac{1}{{\langle\sigma_2\rangle}^b{\langle\sigma_1\rangle}^{b'}}, \ 
\forall(\tau,\xi,\tau_1,\xi_1)\in \mathcal{C}.
\end{displaymath}
Hence, we deduce from \eqref{CEb}, \eqref{algrel2} and \eqref{CEQ}, that 
\begin{equation}\label{psiC}
\|\mathbf 1_{\mathcal{C}}\Psi\|_{L^\infty_{\tau,\xi}(L^2_{\tau_1,\xi_1})}
\lesssim
\sup_{\tau,\xi}{\left[\int\frac{d\xi_1}
{{\langle{(1+\nu\sgn(\xi_1))\xi_1^2-2\xi\xi_1+\sigma}\rangle}^{2\min\{b,b'\}}}\right]}^\frac12
\lesssim1.
\end{equation}
Therefore, we conclude the proof of \eqref{NE2ABC} gathering \eqref{psiA}, \eqref{psiB} and \eqref{psiC}. \\
%
\noindent\textit{Proof of the estimate (\ref{NE2C1})}: By \eqref{xi1mod}, $1\le{|\xi_1|}^2\lesssim|\sigma_1|$ 
in the region ${\mathcal{C}}_1$, thus we have 
\begin{displaymath}
|\Psi(\tau,\xi,\tau_1,\xi_1)|\lesssim 
\frac{{\langle\xi_1\rangle}^{s-s'-\frac12}}
{{\langle\sigma_1\rangle}^{\frac {b'}2}{\langle\sigma_1\rangle}^{\frac {b'}2}}\cdot
\frac{|2\xi_1|^\frac12}{{\langle\sigma_2\rangle}^b{\langle\sigma\rangle}^{1-c}}\lesssim
\frac{|2\xi_1|^\frac12}{{\langle\sigma_1\rangle}^{\frac14}{\langle\sigma_2\rangle}^b{\langle\sigma\rangle}^{1-c}}, 
\ \ \forall(\tau,\xi,\tau_1,\xi_1)\in \mathcal{C}_1, 
\end{displaymath} 
since \eqref{ss'2l} implies $s-s'-1/2<1/2<b'$.
Hence, using \eqref{CE1-b}, \eqref{algrel2} and performing the change of variable 
$\eta:=2\xi_1\xi-(1+\nu\sgn(\xi_1))\xi_1^2+\sigma_1=\sigma-\sigma_2$, we get that 
\begin{equation*}
\|\mathbf 1_{\mathcal{C}_1}\Psi\|_{L^\infty_{\tau_1,\xi_1}(L^2_{\tau,\xi})}
\lesssim
\sup_{\tau_1,\xi_1}{{\langle\sigma_1\rangle}^{-\frac14}
\left[\int_{|\eta|\le2|\sigma_1|}\frac{d\eta}{{\langle{\eta}\rangle}^{2(1-c)}}\right]}^\frac12\\
\lesssim
\sup_{\tau_1,\xi_1}{\langle\sigma_1\rangle}^{-\frac14-\frac12+c}\lesssim 1,
\end{equation*}
since $1/2<c<3/4$. This concludes the proof of \eqref{NE2C1}.\\ 
\noindent\textit{Proof of the estimate (\ref{NE2C2})}: Denoting $\tau:=\tau_1+\tau_2$, $\xi:=\xi_1+\xi_2$ and 
$\sigma, \ \sigma_1, \ \sigma_2$ as before,  we have ${|\xi_1|}^2\lesssim|\sigma_2|$ 
in the region $\widetilde{\mathcal{C}}_2$. Also, $s-s'<1<2b$ by (\ref{ss'2l}). Therefore, 
\begin{displaymath}
|\widetilde{\Psi}(\tau_2,\xi_2,\tau_1,\xi_1)|\lesssim 
\frac{{\langle\xi_1\rangle}^{s-s'}}{{\langle\sigma_2\rangle}^{b}}\cdot
\frac{1}{{\langle\sigma_1\rangle}^{b'}{\langle\sigma\rangle}^{1-c}}\lesssim 
\frac{1}{{\langle\sigma_1\rangle}^{b'}{\langle\sigma\rangle}^{1-c}}, \ 
\forall(\tau_2,\xi_2,\tau_1,\xi_1)\in \widetilde{\mathcal{C}}_2.
\end{displaymath}
Hence, from \eqref{CE1-b}, \eqref{algrel2}, \eqref{CEQ} and $c<3/4$, we conclude that 
\begin{displaymath}
\|\mathbf 1_{\widetilde{\mathcal{C}}_2}\widetilde{\Psi}\|_{L^\infty_{\tau_2,\xi_2}(L^2_{\tau_1,\xi_1})}
\lesssim
\sup_{\tau_2,\xi_2}{\left[\int\frac{1}
{{\langle{(1-\nu\sgn(\xi_1))\xi_1^2+2\xi_1\xi_2+\sigma_2}\rangle}^{2(1-c)}}d\xi_1\right]}^\frac12
\lesssim1.
\end{displaymath}
This finishes the proof of \eqref{NE2} .\hfill$\square$ 
%
%
%
%
\begin{center}
\item\section{Local Well-Posedness}
\end{center}
\setcounter{section}{3} \setcounter{equation}{0} 
Using the new bilinear estimates of the previous section, Theorem \ref{SBO} can be proven, 
with minor adjustments, in the same way that Bekiranov, Ogawa and Ponce proved 
L.W.P. of the system \eqref{SBOsys} for the case $s\ge0,\ s'=s-1/2$. 
In this section, we detail the proof for the convenience of the reader. 
First, we need to state the linear estimates for the Fourier restriction norm method 
(see, e.g., \cite{Gini},  \cite{bop}, \cite{Bo}).\\ 
%
\begin{lem} \label{lem linX}
Let $T\in(0,1), s \in \mathbb R$ and $1/2 < b \le c\le 1$, then the following 
estimates hold:
\begin{eqnarray}
(i)&& \|f\|_{C^0_t(\mathbb R;H^s_x)}\lesssim \|f\|_{X^{s,b}}\ ,\label{lem linX.1}\hspace{10.0cm}\\
\nonumber \\ 
(ii)&& \|\eta(t)e^{it\partial^2_x}\phi\|_{X^{s,b}}\lesssim \|\phi\|_{H^s}\ ,\label{lem linX.2}\\
\nonumber \\ 
(iii)&& \left\|\eta_T(t)\int_0^te^{i(t-t')\partial^2_x}f(t')dt'\right\|_{X^{s,b}}\lesssim 
T^{c-b}\|f\|_{X^{s,c-1}}\ .\label{lem linX.5}
\end{eqnarray}
Similar estimates hold for $e^{-\nu t\mathcal H\partial^2_x}$ and 
$Y^{s,b}$ replacing $e^{it\partial^2_x}$ and $X^{s,b}$, respectively. 
\end{lem} 
%
%
\vspace{0.2cm} 
\noindent \textbf{Proof of Theorem \ref{SBO}.} 
Let $s,s'\in\mathbb R$ satisfy \eqref{ss'1} and \eqref{ss'2}.
Then $-\frac12<\frac{s'-s}2<\frac14$ and we can fix $b,c,b',c'\in\mathbb R$ such that 
\begin{equation*}
\max\left\{1/2\ ,\ (s'-s)/2+1/2\right\}<b<c<\min\left\{3/4\ ,\ (s'-s)/2+1\right\}
\end{equation*}
and 
\begin{equation*}
1/2<b'<c'<\min\left\{3/4-(s'-s)/2\ , \ 3/4\right\}.
\end{equation*}
Thus the hypotheses of Theorems \ref{thm NE1} and \ref{thm NE2} are verified. 
Fix $R>0$, $(\phi,\psi)\in B_R$ and a constant 
$C>0$ 
greater than all the implicit constants in 
the estimates \eqref{NE1}, \eqref{NE2}, 
\eqref{lem linX.2} and \eqref{lem linX.5}, and also greater than $|\alpha|+|\beta|$ . 
Let 
\begin{equation*}
\mathfrak B:=\left\{(u,v)\in X^{s,b}\times Y^{s',b'}\ : \ 
\|(u,v)\|_{\mathfrak B}:=\|u\|_{X^{s,b}}+\|v\|_{Y^{s',b'}}\le 2CR\right\},
\end{equation*} 
which is a complete metric space. For each $T\in(0,1)$ such that 
\begin{equation*}
T^{\min\{c-b\ \!,\ \!c'-b'\}}<(8C^4R)^{-1}, 
\end{equation*}
we consider the map 
$\Xi=\Xi[\phi,\psi,T]:\mathfrak B\to X^{s,b}\times Y^{s',b'}$, $(u,v)\mapsto (\Xi_1(u,v)\ ,\ \Xi_2(u,v))$ 
defined by 
\begin{eqnarray*}
\Xi_1(u,v)&:=&\eta(t)e^{it\partial^2_x}\phi-i\alpha\eta_T(t)\int_0^te^{i(t-t')\partial^2_x}[u(t')v(t')]dt',\\
\Xi_2(u,v)&:=&\eta(t)e^{-\nu t\mathcal H\partial^2_x}\psi+\beta\eta_T(t)
\int_0^te^{-\nu (t-t')\mathcal H\partial^2_x}\partial_x|u(t')|^2dt', 
\end{eqnarray*} 
From the estimates \eqref{NE1}, \eqref{NE2}, \eqref{lem linX.2} and \eqref{lem linX.5}, we conclude that 
\begin{equation*}
\|\Xi(u,v)\|_{X^{s,b}\times Y^{s',b'}}\le CR+(2C^4RT^{\min\{c-b\ \!,\ \!c'-b'\}})(2CR), 
\ \ \ \forall (u,v)\in\mathfrak B, 
\end{equation*}
which means that $\Xi$ maps $\mathfrak B$ on itself, moreover 
\begin{equation*}
\|\Xi(u,v)-\Xi(\tilde u,\tilde v)\|_{\mathfrak B}\le 8C^4RT^{\min\{c-b\ \!,\ \!c'-b'\}}
\|(u,v)-(\tilde u,\tilde v)\|_{\mathfrak B}, 
\ \ \ \ \forall (u,v),(\tilde u,\tilde v)\in\mathfrak B. 
\end{equation*}
Hence, $\Xi:\mathfrak B\to\mathfrak B$ is a contraction and has a unique fixed point. This establishes the 
existence of solution $(u,v)$ satisfying \eqref{uDuhamel} and \eqref{vDuhamel} for every $t\in[-T,T]$, and from 
\eqref{XSobolev} and \eqref{YSobolev} we have 
\begin{equation*}
(u,v)\in C^0([-T,T];H^s(\mathbb R))\times C^0([-T,T];H^{s'}(\mathbb R)). 
\end{equation*}
%
Thus, the flow map data-solution $S$ in \eqref{fluxo} is defined at $(\phi,\psi)\in B_R$ to be the fixed 
point of $\ \Xi[\phi,\psi,T]$. 
From \eqref{lem linX.2}, \eqref{lem linX.5}, \eqref{NE1} and \eqref{NE2}, we get that 
\begin{equation*}
\|S(\phi,\psi)\!-\!S(\tilde\phi,\tilde\psi)\|_\mathfrak B
\le
\lambda \|(\phi,\psi)\!-\!(\tilde\phi,\tilde\psi)\|_{H^s\times H^{s'}}, 
\ \ \ \ \forall(\phi,\psi),(\tilde\phi,\tilde\psi)\!\in \!B_R, 
\end{equation*} 
where $\lambda=C(1-8C^4RT^{\min\{c-b\ \!,\ \!c'-b'\}})^{-1}$. Hence, from \eqref{lem linX.1}, we conclude that 
the flow \eqref{fluxo} is Lipschitz.\\ 
%
\indent Finally, we will prove the uniqueness of the solution in the class $X^{s,b}_T\times Y^{s',b'}_T$. Suppose that 
$(u_1,v_1),(u_2,v_2)\in X^{s,b}\times Y^{s',b'}$ satisfy \eqref{uDuhamel} and \eqref{vDuhamel} for every 
$t\in[-T,T]$.\\ 
\indent Let $T^*\le T$, such that 
\begin{equation}\label{T*} 
2C^3(\|(u_1,v_1)\|_{X^{s,b}\times Y^{s',b'}}+\|(u_2,v_2)\|_{X^{s,b}\times Y^{s',b'}})
{T^*}^{\min\{c-b\ \!,\ \!c'-b'\}}\le\frac12. 
\end{equation}
For any $\epsilon >0$ there exists 
$(\tilde u,\tilde v)\in X^{s,b}\times Y^{s',b'}$ such that 
$\tilde u(t)=u_1(t)-u_2(t)$ and $\tilde v(t)=v_1(t)-v_2(t)$ for every $t\in[-T^*,T^*]$ and 
\begin{equation}\label{utilvtil} 
\|(\tilde u,\tilde v)\|_{X^{s,b}\times Y^{s',b'}}\le\|(u_1,v_1)-(u_2,v_2)\|_{X^{s,b}_{T^*}\times Y^{s',b'}_{T^*}}+\epsilon. 
\end{equation} 
Therefore, for every $t\in[-T^*,T^*]$, 
\begin{eqnarray*}
u_1(t)-u_2(t)&=&-i\alpha\eta_{T^*}(t)\int_0^te^{i(t-t')\partial^2_x}[\tilde u(t')v_1(t')+u_2(t')\tilde v(t')]dt',\\
v_1(t)-v_2(t)&=&\beta\eta_{T^*}(t)\int_0^te^{-\nu (t-t')\mathcal H\partial^2_x}\partial_x
[\tilde u(t')\overline{u_1(t')}+u_2(t')\overline{\tilde u(t')}]dt'. 
\end{eqnarray*} 
Thus from \eqref{lem linX.5}, \eqref{NE1} and \eqref{NE2} yields 
\begin{eqnarray}
\|u_1-u_2\|_{X^{s,b}_{T^*}}
&\le&
\left\|-i\alpha\eta_{T^*}(t)\int_0^te^{i(t-t')\partial^2_x}[\tilde u(t')v_1(t')+u_2(t')\tilde v(t')]dt'\right\|_{X^{s,b}}
\nonumber\\ 
&\le&
C^3{T^*}^{c-b}\left(\|\tilde u\|_{X^{s,b}}\|v_1\|_{Y^{s',b'}}+\|u_2\|_{X^{s,b}}\|\tilde v\|_{Y^{s',b'}}\right)\label{u1u2} 
\end{eqnarray} 
and 
\begin{equation}\label{v1v2}
\|v_1-v_2\|_{Y^{s',b'}_{T^*}}
\le 
C^3{T^*}^{c'-b'}\left(\|\tilde u\|_{X^{s,b}}\|u_1\|_{X^{s,b}}+\|u_2\|_{X^{s,b}}\|\tilde u\|_{X^{s,b}}\right). 
\end{equation}
Combining \eqref{T*}, \eqref{u1u2} and \eqref{v1v2} we have 
\begin{equation}\label{u1util}
\|(u_1,v_1)-(u_2,v_2)\|_{X^{s,b}_{T^*}\times Y^{s',b'}_{T^*}}
\le 
\frac12\|(\tilde u,\tilde v)\|_{X^{s,b}\times Y^{s',b'}}. 
\end{equation} 
From \eqref{utilvtil} and \eqref{u1util}, we conclude that 
$\|(u_1,v_1)-(u_2,v_2)\|_{X^{s,b}_{T^*}\times Y^{s',b'}_{T^*}}\le\epsilon$. Hence, since $\epsilon$ is arbitrary, 
$(u_1,v_1)=(u_2,v_2)$ on $[-T^*,T^*]$. Using translations in time, one can repeat this argument a finite number 
of times to conclude that $(u_1,v_1)=(u_2,v_2)$ on $[-T,T]$. \hfill $\square$ 
%
%
%
%
\begin{center}
\item\section{Ill-Posedness Results}
\end{center}
\setcounter{section}{4} \setcounter{equation}{0} 
Suppose that there exists $T\!\!>\!\!0$ such that the Cauchy problem 
\eqref{SBOsys} is locally well-posed in $H^s(\mathbb R)\!\times\!H^{s'}(\mathbb R)$, 
in the time interval $[-T,T]$. 
Suppose also that there exists $t\!\in\![-T,0)\cup(0,T]$ such that the associated flow map data-solution 
\eqref{fluxot} is 
two times Fréchet differentiable 
at zero. 
Then the second Fréchet derivative of $S^t$ at zero belongs to 
$\mathcal B$, 
the normed space of bounded bilinear applications from $(H^s\!\times\!H^{s'})^2$ to $H^s\!\times\!H^{s'}$. 
In particular, we have the following estimate for the second Gâteaux derivative of $S^t$ at zero, 
\begin{equation}\label{eqvC2ill}
\left\|\frac{\partial^2S^t}{\partial (\phi,\psi)^2}(0,0)\right\|_{H^s\!\times\!H^{s'}}
\le 
\|D^2S^t(0,0)\|_{\mathcal B}\cdot\|(\phi,\psi)\|^2_{H^s\!\times\!H^{s'}} 
\ , \ \ \forall \phi,\psi\in\mathcal S(\mathbb R).
\end{equation}
We will denote $(u_{\phi,\psi}(t),v_{\phi,\psi}(t)):=S^t(\phi,\psi)$. 
This means that $(u_{\phi,\psi}(t),v_{\phi,\psi}(t))$ is a solution of the associated integral equations
\begin{eqnarray}
u_{\phi,\psi}(t)&=&e^{it\partial^2_x}\phi-i\alpha\int^t_0
e^{i(t-t')\partial^2_x}\left(u_{\phi,\psi}(t')\cdot v_{\phi,\psi}(t')\right)dt',\label{uweqint}\\
v_{\phi,\psi}(t)&=&e^{-\nu t\mathcal H\partial^2_x}\psi+\beta\int^t_0
e^{-\nu(t-t')\mathcal H\partial^2_x}(\partial_x|u_{\phi,\psi}(t')|^2)dt'\label{vweqint}.
\end{eqnarray} 
Since $(u_{0,0}(t) , v_{0,0}(t))=S^t(0,0)=(0,0)$, we have 
\begin{equation*}
\frac{\partial S^t}{\partial(\phi,\psi)}(0,0)=
\left(\frac{\partial u_{0,0}}{\partial(\phi,\psi)}(t)\ ,\ \frac{\partial v_{0,0}}{\partial(\phi,\psi)}(t)\right)
=
\left(e^{it\partial^2_x}\phi\ ,\ e^{-\nu t\mathcal H\partial^2_x}\psi\right), 
\ \ \ \forall\phi,\psi\in\mathcal S(\mathbb R).
\end{equation*}
Thus, using \eqref{uweqint} to compute the second Gâteaux derivative of $u$ at zero, in direction 
$(\phi,\psi)\in\mathcal S(\mathbb R)\times\mathcal S(\mathbb R)$, yields 
\begin{equation*}
\frac{\partial^2u_{0,0}}{{\partial (\phi,\psi)}^2}(t)=
-2i\alpha\int^t_0e^{i(t-t')\partial^2_x}\left(e^{it'\partial^2_x}\phi\cdot 
e^{-\nu t'\mathcal H\partial^2_x}\psi\right)dt'.
\end{equation*}
Therefore, denoting  $\xi_2:=\xi-\xi_1$, we have 
\begin{eqnarray*}
\left\|\frac{\partial^2u_{0,0}}{{\partial (\phi,\psi)}^2}(t)\right\|_{H^s}
&=&
\left\|2\alpha{\langle\xi\rangle}^s\int^t_0e^{-i(t-t')\xi^2}
\left((e^{it'\partial^2_x}\phi)\ \widehat{}\ \ast (e^{-\nu t'\mathcal H\partial^2_x}\psi)\ \widehat{}\ \right)(\xi)
dt'\right\|_{L^2_\xi}\\
&=&
\left\|2\alpha{\langle\xi\rangle}^s\int^t_0e^{it'\xi^2}\int e^{-it'(\xi_2^2+\nu|\xi_1|\xi_1)}\widehat{\phi}(\xi_2)
\widehat{\psi}(\xi_1)d\xi_1dt'\right\|_{L^2_\xi}\\
&=&
\left\|\int_0^t\int\Theta(t',\xi,\xi_1)f(\xi_2)g(\xi_1)d\xi_1dt\right\|_{L^2_\xi},
\end{eqnarray*}
where 
$f(\xi_2):={\langle\xi_2\rangle}^s\widehat{\phi}(\xi_2)$, 
$g(\xi_1):={\langle\xi_1\rangle}^{s'}\widehat{\psi}(\xi_1)$ and 
\begin{equation}\label{Thetadef}
\Theta(t',\xi,\xi_1):=\frac{2|\alpha|{\langle\xi\rangle}^s}{{\langle\xi_2\rangle}^{s}{\langle\xi_1\rangle}^{s'}}
\cdot e^{it'(\xi^2-\xi_2^2-\nu|\xi_1|\xi_1)}\ . 
\end{equation}
Hence, the assumption that the flow map \eqref{fluxot} is $C^2$ at zero implies  
\begin{equation}\label{eqv2C2ill}
\left\|\int^t_0\!\!\int\!\Theta(t'\!,\xi,\xi_1)f(\xi_2)g(\xi_1)d\xi_1dt'\right\|_{L^2_\xi}
\!\!\le
\!\|D^2S^t(0,0)\|_{\mathcal B} 
(\|f\|_{L^2}\!+\!\|g\|_{L^2})^2, 
\forall f,g\!\in\!\mathcal S(\mathbb R). 
\end{equation}
Similarly, differentiating the equation 
\eqref{vweqint} twice, in direction $(\phi,0)\in\mathcal S(\mathbb R)\times\mathcal S(\mathbb R)$, yields
\begin{equation*}
\frac{\partial^2v_{0,0}}{\partial (\phi,0)^2}(t)=
2\beta\int^t_0e^{-\nu(t-t')\mathcal H\partial^2_x}\partial_x\left(e^{it'\partial^2_x}\phi
\cdot\overline{e^{it'\partial^2_x}\phi}\right)dt'. 
\end{equation*}
Thus, that assumption for the flow map \eqref{fluxot} also implies 
\begin{equation}\label{eqvC2ill2}
\left\|\int^t_0\int\Upsilon(t',\xi,\xi_1)f(\xi_2)\overline{f(-\xi_1)}d\xi_1dt'\right\|_{L^2_\xi}
\le
\|D^2S^t(0,0)\|_{\mathcal B}\cdot\|f\|^2_{L^2}, \ \ \ \forall f\in\mathcal S(\mathbb R),
\end{equation}
where, 
\begin{equation}\label{Upsilondef} 
\Upsilon(t',\xi,\xi_1):=\frac{2|\beta||i\xi|{\langle\xi\rangle}^{s'}}{{\langle\xi_2\rangle}^s{\langle\xi_1\rangle}^s}
\cdot e^{it'(\nu|\xi|\xi-\xi_2^2+\xi_1^2)}. 
\vspace{0.3cm} 
\end{equation}
Next, we will state an elementary result that will be useful in the proofs of Theorems \ref{C2ill}, \ref{C2ill2}, 
\ref{thm notNE} and \ref{thm endpoint}. 
\vspace{0.2cm} 
\begin{lem}\label{lemconvbaixo}
Let $A, B, R \subset \mathbb R^n$. 
If $R - B \subset A$ then
\begin{equation}\label{convbaixo}
\|\mathbf 1_R\|_{L^2(\mathbb R^n)} \|\mathbf 1_B\|_{L^1(\mathbb R^n)}
\le
 \|\mathbf 1_A \ast \mathbf 1_B\|_{L^2(\mathbb R^n)}.
\end{equation}
\end{lem}
\noindent \textbf{Proof.} 
If $R - B \subset A$, then 
$$\mathbf 1_A \ast \mathbf 1_B(x)
=\int_A \mathbf 1_B(x-y)dy
=\int_A \mathbf 1_{x-B}(y)dy
\ge \mathbf 1_R(x)\|\mathbf 1_B\|_{L^1(\mathbb R^n)}, \ \ \forall x\in\mathbb R^n,$$ 
taking the $L^2$-norm, (\ref{convbaixo}) follows. 
\hfill $\square$ 
%
\vspace{0.5cm} \\ 
%
\noindent \textbf{Proof of Theorem \ref{C2ill}.} 
\textit{(i)} It is enough to show that \eqref{eqv2C2ill} or \eqref{eqvC2ill2} fails.\\ 
\textit{Case $s'<-1/2$}: In this case, \eqref{eqv2C2ill} fails. 
Indeed, for each $N\in \mathbb N$, define 
\begin{eqnarray*}
A_N&:=& \{\xi_1\in\mathbb R \ : \ |(1+|\nu|)\xi_1+\sgn(\nu)N|<
(1+|\nu|)(4\langle t\rangle N)^{-1}\},\\
B_N&:=& \{\xi_2\in\mathbb R \ : \  |2\xi_2-\sgn(\nu)N|< (4\langle t\rangle N)^{-1}\}.
\end{eqnarray*} 
For $N$ sufficiently large (precisely $N>1+|\nu|$), 
we have that 
\begin{equation}\label{freqN}
\langle\xi_1\rangle\sim\langle\xi_2\rangle\sim\langle\xi_1+\xi_2\rangle\sim N , 
\ \ \ \ \ \ \forall \xi_1\in A_N,\ \forall \xi_2\in B_N, 
\end{equation} 
since $1+|\nu|\ne2$. Moreover, $\sgn(\xi_1)=-\sgn(\nu)$ for all $\xi_1\in A_N$. Thus, we also have 
\begin{equation}\label{ressont1}
|(\xi_1+\xi_2)^2-\nu|\xi_1|\xi_1-\xi_2^2|=|\xi_1|\cdot|(1+|\nu|)\xi_1+2\xi_2|
<\frac{2N}{1+|\nu|}\cdot\frac{2+|\nu|}{4\langle t\rangle N}
\le\frac1{|t|}. 
\end{equation}
Observe that $\cos(x)\ge 1/2$ for $|x|\le 1$. 
Hence, we deduce from \eqref{Thetadef}, \eqref{freqN} and \eqref{ressont1} that 
\begin{equation}\label{ReTheta}
Re\left(\Theta(t',\xi_1+\xi_2,\xi_1)\right)\gtrsim \frac1{N^{s'}}\ , 
\ \ \ \ \ \forall|t'|\le|t|, \forall\xi_1\in A_N, \forall\xi_2\in B_N. 
\end{equation}
Now, taking $f_N,g_N\in\mathcal S(\mathbb R)$ such that $\mathbf 1_{A_N}\le g_N$, $\mathbf 1_{B_N}\le f_N$, 
$\|g_N\|_{L^2}\le2\|\mathbf 1_{A_N}\|_{L^2}$ and $\|f_N\|_{L^2}\le2\|\mathbf 1_{B_N}\|_{L^2}$, 
and using \eqref{ReTheta}, we get that 
\begin{eqnarray}
\left|\int^t_0\!\!\int\!\!\Theta(t',\xi,\xi_1)f_N(\xi\!-\!\xi_1)g_N(\xi_1)d\xi_1dt'\right|
\!\!\!&\ge&\!\!\!
\left|\int^t_0\!\!\int\!\!Re(\Theta(t',\xi,\xi_1)) 
\mathbf 1_{B_N}(\xi\!-\!\xi_1)\mathbf 1_{A_N}(\xi_1)d\xi_1dt'\right|\nonumber\\ 
&\gtrsim&\!\!\!
|t|\cdot\frac{\mathbf 1_{A_N}\ast\mathbf 1_{B_N}(\xi)}{N^{s'}}\ , 
\ \ \ \ \ \ \ \ \ \ \ \ \ \ \ \forall\xi\in\mathbb R.\label{estTheta}
\end{eqnarray}
Combining \eqref{estTheta} with \eqref{eqv2C2ill}, yields 
\begin{equation}\label{absur1C2}
|t|\cdot\frac{\|\mathbf 1_{A_N}\ast\mathbf 1_{B_N}\|_{L^2}}{N^{s'}}
\lesssim
\|D^2S^t(0,0)\|_{\mathcal B}\cdot 
(\|\mathbf 1_{A_N}\|_{L^2}+\|\mathbf 1_{B_N}\|_{L^2})^2.
\end{equation}
On the other hand, defining
\begin{equation*}
R_N\ :=\ \{\xi\in\mathbb R \ : \ |\xi+b_{\nu}\sgn(\nu)N|<{(8\langle t\rangle N)}^{-1}\},
\end{equation*}
where $b_{\nu}=\frac{1}{1+|\nu|}-\frac12\ne0$, we have 
$R_N - B_N \subset A_N$. 
Hence, from \eqref{convbaixo} and \eqref{absur1C2}, we conclude that 
\begin{equation}\label{absur2C2}
|t|\cdot\frac{N^{-\frac12}N^{-1}}{N^{s'}}\lesssim
\frac{\|D^2S^t(0,0)\|_{\mathcal B}}{N},
\end{equation}
which is false in the case $s'<-1/2$, since $N$ can be chosen arbitrarily large. 
\vspace{0.2cm} \\ 
\noindent\textit{Case $s'>2s-1/2$}: 
For this case, we will show that \eqref{eqvC2ill2} fails, using the 
same ideas used in the previous case. 
For $N\in \mathbb N$ sufficiently large (precisely $N>|1-|\nu||^{-1}$), define 
\begin{eqnarray*}
A_N&:=& \{\xi_1\in\mathbb R \ : \ |a_\nu\xi_1+\sgn(\nu)(1+|\nu|)N|<(c_tN)^{-1}\},\\
B_N&:=& \{\xi_2\in\mathbb R \ : \ |a_\nu\xi_2+\sgn(\nu)(1-|\nu|)N|<(2c_tN)^{-1}\},\\
R_N&:=& \{\xi\in\mathbb R \ : \  |a_\nu\xi+2\sgn(\nu)N|<(2c_tN)^{-1}\},
\end{eqnarray*}
where $a_\nu:=|1-|\nu||\cdot|1+|\nu||\ne0$ and $c_t:=1+8|t|(1-|\nu|)^{-2}$. Then $R_N - B_N \subset A_N$.\\ 
And if $\xi_1\!\in\! A_N$ and $\xi_2\!\in\! B_N$, then 
$\langle\xi_1\rangle\!\sim\!\langle\xi_2\rangle\!\sim\!\langle\xi_1+\xi_2\rangle\!\sim\! N$, 
$\sgn(\xi_1+\xi_2)\!=\!-\!\sgn(\nu)$ and 
\begin{eqnarray}
\left|\nu(\xi_1+\xi_2)|\xi_1+\xi_2|-\xi_2^2+\xi_1^2\right|
&=&
|\xi_1+\xi_2|\cdot|(1-|\nu|)\xi_1-(1+|\nu|)\xi_2|\nonumber\\
&<&\frac{4N}{a_\nu}\cdot\frac{2}{c_t|1-|\nu||N}
\le\frac1{|t|}\ .\label{ressont2} 
\end{eqnarray} 
Following the arguments used in \eqref{freqN}-\eqref{ReTheta}, 
we get from \eqref{Upsilondef} and \eqref{ressont2} that 
\begin{equation}\label{ReUpsilon}
Re\left(\Upsilon(t',\xi_1+\xi_2,\xi_1)\right)\gtrsim \frac{N^{s'+1}}{N^{2s}}\ , 
\ \ \ \ \forall|t'|\le|t|, \forall\xi_1\in A_N, \forall\xi_2\in B_N. 
\end{equation}
Now, taking $f_N\!\in\!\mathcal S(\mathbb R)$ such that $\mathbf 1_{-A_N\cup B_N}\le f_N$ and 
$\|f_N\|_{L^2}\le2\|\mathbf 1_{-A_N\cup B_N}\|_{L^2}\lesssim 
N^{-\frac12}$, 
yields 
\begin{equation*} 
f_N(\xi-\xi_1)\overline{f_N(-\xi_1)}
\ge 
\mathbf 1_{-A_N\cup B_N}(\xi-\xi_1)\mathbf 1_{-A_N\cup B_N}(-\xi_1) 
\ge 
\mathbf 1_{B_N}(\xi-\xi_1)\mathbf 1_{A_N}(\xi_1), 
\end{equation*} 
for all $\xi,\xi_1\in\mathbb R$. 
Thus, similarly to \eqref{estTheta}, we deduce from \eqref{ReUpsilon} that   
\begin{equation}\label{estUpsilon}
\left|\int^t_0\!\int\Upsilon(t',\xi,\xi_1)f_N(\xi-\xi_1)\overline{f_N(-\xi_1)}d\xi_1dt'\right|
\gtrsim
|t|\cdot\frac{\mathbf 1_{A_N}\ast\mathbf 1_{B_N}(\xi)\cdot N^{s'+1}}{N^{2s}}, 
\ \ \ \ \forall\xi\in\mathbb R.
\end{equation}
Combining \eqref{estUpsilon}, \eqref{eqvC2ill2} and \eqref{convbaixo}, we conclude 
\begin{equation}\label{absur4C2}
|t|\cdot\frac{N^{s'+1}\cdot N^{-\frac12}\cdot N^{-1}}{N^{2s}}
\lesssim
\frac{\|D^2S^t(0,0)\|_{\mathcal B}}{N},
\end{equation}
which is false in the case $2s-1/2<s'$, since $N$ can be chosen arbitrarily large. 
\vspace{0.2cm}\\ 
\noindent\textit{(ii)} If the map \eqref{fluxo} is $C^2$ at zero then \eqref{eqv2C2ill} 
and \eqref{eqvC2ill2} hold for every $t\in[-T,T]$ and 
\begin{equation}\label{fluxnC2}
\sup_{t\in[-T,T]}\|D^2S^t(0,0)\|_{\mathcal B}
<\infty.
\end{equation}
Thus, it is enough to show that \eqref{fluxnC2} fails for $\left|s'\!-\!(s\!-\!1/2)\right|\!>\!3/2$, i.e., 
for $s'<s-2$ or $s+1<s'$. 
Indeed, for each $N\in\mathbb N$, defining 
\begin{eqnarray*}
A_N&:=& \{\xi_1\in\mathbb R \ : \ |\xi_1-N|<1/2\},\\
B_N&:=& \{\xi_2\in\mathbb R \ : \  |\xi_2|<1/4\},\\
R_N&:=& \{\xi\in\mathbb R \ : \  |\xi-N|<1/4\}, 
\end{eqnarray*}
we have $R_N - B_N \subset A_N$. Also, if $\xi_1\in A_N$ and $\xi_2\in B_N$ then  
\begin{displaymath}
\langle\xi_1\rangle\sim N, \ \ \langle\xi_2\rangle\sim1, \ \ \langle\xi_1+\xi_2\rangle\sim N, 
\end{displaymath}
and 
\begin{equation*}
\left|(\xi_1+\xi_2)^2-\nu|\xi_1|\xi_1-\xi_2^2\right|=
|\xi_1|\cdot|(1-\nu\sgn(\xi_1))\xi_1+2\xi_2|
<6(1+|\nu|)N^2. 
\end{equation*}
In addition, for $N>(6(1+|\nu|)T)^{-\frac12}$, we define 
$t_N:=(6(1+|\nu|)N^2)^{-1}\!\in\!(0,T]$. 
Therefore, following the arguments used in \eqref{freqN}-\eqref{absur2C2}, we get that 
\begin{equation*}
N^{s-2-s'}\lesssim 
\frac{N^s\cdot t_N}{N^{s'}}
\lesssim
\|D^2S^{t_N}(0,0)\|_{\mathcal B}, 
\end{equation*} 
contradicting \eqref{fluxnC2} when $s'<s-2$ (since $N$ can be chosen arbitrarily large).\\ 
Moreover, for $\xi_1\in A_N$ and $\xi_2\in B_N$, we have 
\begin{equation*}
\left|\nu(\xi_1+\xi_2)|\xi_1+\xi_2|-\xi_2^2+\xi_1^2\right|
<
6(1+|\nu|)N^2. 
\end{equation*}
Now following \eqref{ressont2}-\eqref{absur4C2} we conclude that 
\begin{equation*}
N^{s'-s-1}
\lesssim 
\frac{N^{s'+1}\cdot t_N}{N^s}
\lesssim
\|D^2S^{t_N}(0,0)\|_{\mathcal B}, 
\end{equation*} 
contradicting \eqref{fluxnC2} when $s+1<s'$. 
This finishes the proof of the theorem. 
\hfill$\square$ 
\vspace{0.5cm} \\ 
%
\noindent \textbf{Proof of Theorem \ref{C2ill2}.} The proof use the same arguments used in the proof of 
Theorem \ref{C2ill}\textit{(i)}. Suppose that we have some $t\in [-T,0)\cup(0,T]$ such that the flow 
map \eqref{fluxot} is $C^2$ at zero. 
For $N\in\mathbb N$, defining  
\begin{eqnarray*}
A_N&:=& \{\xi_1\in\mathbb R \ : \ |\xi_1-\sgn(\nu)N|<(2\langle t\rangle N)^{-1}\},\\
B_N&:=& \{\xi_2\in\mathbb R \ : \  |\xi_2|<(4\langle t\rangle N)^{-1}\},\\
R_N&:=& \{\xi\in\mathbb R \ : \  |\xi-\sgn(\nu)N|<(4\langle t\rangle N)^{-1}\}, 
\end{eqnarray*}
we have $R_N - B_N \subset A_N$. Also, if $\xi_1\in A_N$, $\xi_2\in B_N$ then 
$\langle\xi_1\rangle\sim N$, $\langle\xi_2\rangle\sim1$, $\langle\xi_1+\xi_2\rangle\sim N$, 
$\sgn(\xi_1)=\sgn(\nu)$ and 
\begin{equation*}
\left|(\xi_1+\xi_2)^2-\nu|\xi_1|\xi_1-\xi_2^2\right|=|2\xi_1\xi_2|
<4N\cdot(4\langle t\rangle N)^{-1}\le|t|^{-1}. 
\end{equation*}
Following the arguments used in \eqref{freqN}-\eqref{absur2C2}, we deduce from \eqref{eqv2C2ill} 
that 
\begin{displaymath}
|t|\cdot N^{s-\frac12-s'}\lesssim\|D^2S^t(0,0)\|_{\mathcal B},\ \ \ \ \forall N\in\mathbb N. 
\end{displaymath}
Hence $s'\ge s-1/2$. 
On the other hand, defining  
\begin{eqnarray*}
A_N&:=& \{\xi_1\in\mathbb R \ : \ |\xi_1+\sgn(\nu)N|<(2\langle t\rangle N)^{-1}\},\\
B_N&:=& \{\xi_2\in\mathbb R \ : \  |\xi_2|<(4\langle t\rangle N)^{-1}\},\\
R_N&:=& \{\xi\in\mathbb R \ : \  |\xi+\sgn(\nu)N|<(4\langle t\rangle N)^{-1}\},
\end{eqnarray*}
we have $R_N - B_N \subset A_N$. 
Also, if $\xi_1\in A_N$, $\xi_2\in B_N$ then 
$\langle\xi_1\rangle\sim N$, $\langle\xi_2\rangle\sim1$, $\langle\xi_1+\xi_2\rangle\sim N$, 
$\sgn(\xi_1+\xi_2)=-\sgn(\nu)$ and 
\begin{equation*}
\left|\nu(\xi_1+\xi_2)|\xi_1+\xi_2|-\xi_2^2+\xi_1^2\right|=|2(\xi_1+\xi_2)\xi_2|
<4N\cdot(4\langle t\rangle N)^{-1}\le|t|^{-1}. 
\end{equation*}
Now following the arguments used in \eqref{ressont2}-\eqref{absur4C2}, we get from \eqref{eqvC2ill2} 
that 
\begin{displaymath}
|t|\cdot N^{s'-s+\frac12}\lesssim\|D^2S^t(0,0)\|_{\mathcal B},\ \ \ \ \forall N\in\mathbb N. 
\end{displaymath} 
Hence $s'\le s-1/2\le s'$. \\ 
Finally, we will conclude that $s\ge0$. Defining, 
\begin{eqnarray*}
A_N&:=& \{\xi_1\in\mathbb R \ : \ |\xi_1+\sgn(\nu)N|<(8\langle t\rangle N)^{-1}\},\\
B_N&:=& \{\xi_2\in\mathbb R \ : \  |\xi_2-\sgn(\nu)N|< (16\langle t\rangle N)^{-1}\},\\ 
R_N&:=& \{\xi\in\mathbb R \ : \ |\xi|<{(16\langle t\rangle N)}^{-1}\},
\end{eqnarray*}
we have $R_N - B_N \subset A_N$. Also, if $\xi_1\in A_N$, $\xi_2\in B_N$ then 
$\langle\xi_1\rangle\sim\langle\xi_2\rangle\sim N$, $\langle\xi_1+\xi_2\rangle\sim 1$, $\sgn(\xi_1)=-\sgn(\nu)$ 
and 
\begin{equation*}
\left|(\xi_1+\xi_2)^2-\nu|\xi_1|\xi_1-\xi_2^2\right|=|2\xi_1(\xi_1+\xi_2)|
<4N\cdot(4\langle t\rangle N)^{-1}\le|t|^{-1}. 
\end{equation*}
Thus, similarly to \eqref{absur2C2}, we conclude from \eqref{eqv2C2ill} that 
$|t|\cdot N^{-\frac12-s'-s}\lesssim\|D^2S^t(0,0)\|_{\mathcal B}$, for every $N\in\mathbb N$. 
Hence $-2s=-1/2-s'-s\le0$, and this finishes the proof.\hfill $\square$ 
\vspace{0.2cm} \\ 
%
\indent 
We finish this section giving some results about the remaining regions. 
For the \emph{non-resonant} case, Theorem \ref{thm notNE} states that, 
in a part of the remaining region, the L.W.P. of \eqref{SBOsys} can not be obtained by using 
the method of proof employed in this paper. 
Note that, in the case where $\nu=0$, the method fails in the whole remaining region. 
In the \emph{resonant} case, Theorem \ref{thm endpoint} ensures that the method used in \cite{per} can 
not provide L.W.P. for \eqref{SBOsys} at the end-point. \\ 
\indent 
Our proofs of Theorems \ref{thm notNE} and \ref{thm endpoint} follow the arguments used 
by Kenig, Ponce and Vega in \cite{KPV} 
to prove that their $X^{s,b}$ bilinear estimate for KdV equation fails for $s\!<\!-3/4$. 
But in our setting, Lemma \ref{lemconvbaixo} allows to give slightly more direct proofs. 
%
\begin{thm}\label{thm notNE} 
Let $|\nu|\ne1$ and $s,s',c,c'\in\mathbb R$. For every $c',c>1/2$,
\begin{itemize}
\item[(i)] 
the bilinear estimate \eqref{NE1} fails for $s+1/2\le s'$;
\item[(ii)] 
the bilinear estimate \eqref{NE2} fails for $s'\le s-3/2$;
\item[(iii)] 
the bilinear estimate \eqref{NE2} fails for $s'\le s-1$, when $\nu=0$.
\end{itemize}
\end{thm}
%
\noindent \textbf{Proof.} 
\textit{(i)} Recalling the notations of the proof of Theorem 
\ref{thm NE1}, we just have to show that \eqref{NE1eqv} fails when $s+1/2\le s'$. 
For $N\in\mathbb N$, defining 
\begin{eqnarray*}
A_N&:=& \{(\tau_1,\xi_1)\in\mathbb R^2 \ : \ |\xi_1-N|<N^{-1}, \ |\sigma_1|< 6 \},\\
B_N&:=& \{(\tau_2,\xi_2)\in\mathbb R^2 \ : \ |\xi_2|< (2N)^{-1}, \ |\sigma_2|<1 \},\\ 
R_N&:=&\{(\tau,\xi)\in\mathbb R^2 \ : \ |\xi-N|<{(2N)}^{-1}, \ |\tau+\xi^2-2\xi N|< 1 \}, 
\end{eqnarray*}
we have $R_N - B_N \subset A_N$, since $\sigma_1+\sigma_2=\tau+\xi^2-2\xi\xi_1$. 
Moreover, for all $(\tau_1,\xi_1)\in A_N$ and $(\tau_2,\xi_2)\in B_N$, 
\begin{equation*}
\langle\xi_1\rangle\sim N,\ \ \ \ \langle\xi_2\rangle\sim 1,\ \ \ \ \langle\xi\rangle\sim N,\ \ \ \ 
\langle\sigma\rangle\lesssim N^2.
\end{equation*}
Therefore, for all $(\tau,\xi,\tau_1,\xi_1)\in \mathbb R^4$,  
\begin{equation}\label{phiNE1des}
\frac{N^{s'+1}\cdot\mathbf 1_{B_N}(\tau_2,\xi_2)\mathbf 1_{A_N}(\tau_1,\xi_1)}{N^{2(1-c')}\cdot N^s}\lesssim
|\Phi(\tau,\xi,\tau_1,\xi_1)\mathbf 1_{B_N}(\tau_2,\xi_2)\mathbf 1_{A_N}(\tau_1,\xi_1)|.
\end{equation} 
Now, taking $f_N,g_N\!\in\!\mathcal S(\mathbb R^2)$ such that $\mathbf 1_{A_N}\le g_N$, 
$\mathbf 1_{B_N}\le f_N$, 
$\|g_N\|_{L^2}\lesssim\|\mathbf 1_{A_N}\|_{L^2}$, $\|f_N\|_{L^2}\lesssim\|\mathbf 1_{B_N}\|_{L^2}$
and combining \eqref{NE1eqv}, \eqref{phiNE1des} and \eqref{convbaixo}, yields the estimate
\begin{equation*}
\frac{N^{s'+1}\cdot N^{-\frac12}\cdot N^{-1}}{N^{2(1-c')}\cdot N^s}\lesssim N^{-\frac12}\cdot N^{-\frac12}, 
\end{equation*}
which is false for $N$ sufficiently large whenever $s+1/2\le s'$ and $c'>1/2$. 
\vspace{0.2cm}\\ 
\textit{(ii)} Recalling the notations of the proof of Theorem \ref{thm NE2}, we just have to show that 
\eqref{NE2eqv} fails when $s'\le s-3/2$. For $N\in\mathbb N$, defining 
\begin{eqnarray*}
A_N&:=& \{(\tau_1,\xi_1)\in\mathbb R^2 \ : \ |\xi_1-\sgn(\nu)N|<N^{-1}, \ |\sigma_1|< 7(1+|\nu|)\},\\
B_N&:=& \{(\tau_2,\xi_2)\in\mathbb R^2 \ : \ |\xi_2|< (2N)^{-1}, \ |\sigma_2|<1\},\\
R_N&:=& \{(\tau,\xi)\in\mathbb R^2 \ : \ |\xi-\sgn(\nu)N|< (2N)^{-1}, \ |\tau+\xi^2+a_\nu\sgn(\nu)N\xi|<1\},
\end{eqnarray*}
where $a_\nu:=|\nu|-1$, we have $R_N - B_N \subset A_N$, since 
\begin{equation*} 
\sigma_1+\sigma_2=[\tau+\xi^2+a_\nu\sgn(\nu)N\xi]+[(1+|\nu|)\xi_1(\xi_1-\xi)]+[a_\nu(\xi_1-\sgn(\nu)N)\xi].
\end{equation*} 
Arguing as in the previous case, we get from \eqref{NE2eqv} the following estimate 
\begin{displaymath}
\frac{N^s\cdot N^{-\frac12}\cdot N^{-1}}{N^{2(1-c)}\cdot N^{s'}}\lesssim N^{-\frac12}\cdot N^{-\frac12}, 
\end{displaymath}
which is false for $N$ sufficiently large when $s'\le s-3/2$ and $c>1/2$. 
\vspace{0.2cm}\\ 
\textit{(iii)} Recalling the notations of the proof of Theorem \ref{thm NE2} and defining for each $N\in\mathbb N$,  
\begin{eqnarray*}
A_N&:=& \{(\tau_1,\xi_1)\in\mathbb R^2 \ : \ |\xi_1-N|<1, \ |\sigma_1|< 3 \},\\
B_N&:=& \{(\tau_2,\xi_2)\in\mathbb R^2 \ : \  |\xi_2|< 1/2, \ |\sigma_2|<1 \},\\
R_N&:=& \{(\tau,\xi)\in\mathbb R^2 \ : \  |\xi-N|< 1/2, \ |\tau|<1 \},
\end{eqnarray*}
we have $R_N - B_N \subset A_N$, since $\sigma_1+\sigma_2=\tau+\xi_2^2$ in the particular case $\nu=0$. 
Arguing as in the case $(i)$, we get from \eqref{NE2eqv} the following estimate 
\begin{displaymath}
\frac{N^s}{N^{2(1-c)}\cdot N^{s'}}\lesssim 1, 
\end{displaymath}
which is false for $N$ sufficiently large when $s'\le s-1$ and $c>1/2$. \hfill$\square$
%
%
\begin{thm}\label{thm endpoint}
Let $|\nu|=1$, $(s,s')=(0-1/2)$ and $c,b,b'\in\mathbb R$. The estimate \eqref{NE2} fails for every $c>1/2$. 
\end{thm}
\noindent \textbf{Proof.} 
Recalling the notations of the proof of Theorem \ref{thm NE2}, we just have to show that 
\eqref{NE2eqv} fails. For $N\in\mathbb N$, defining 
\begin{eqnarray*}
A_N&:=& \{(\tau_1,\xi_1)\in\mathbb R^2 \ : \ |\xi_1+\sgn(\nu)N|<1/2, \ |\sigma_1|< 1\},\\
B_N&:=& \{(\tau_2,\xi_2)\in\mathbb R^2 \ : \ |\xi_2-\sgn(\nu)N|< 1/4, \ |\sigma_2|<1/3\},\\
R_N&:=& \{(\tau,\xi)\in\mathbb R^2 \ : \ |\xi|< 1/4, \ |\sigma+2\sgn(\nu)N\xi|<1/3\},
\end{eqnarray*} 
we have $R_N - B_N \subset A_N$. 
Moreover, for all $(\tau_1,\xi_1)\in A_N$ and $(\tau_2,\xi_2)\in B_N$, 
\begin{equation*}
\langle\xi_1\rangle\sim N,\ \ \ \ \langle\xi_2\rangle\sim N,\ \ \ \ \langle\xi\rangle\sim 1,\ \ \ \ 
\langle\sigma\rangle\lesssim N.
\end{equation*}
Therefore, for all $(\tau,\xi,\tau_1,\xi_1)\in \mathbb R^4$,  
\begin{equation}\label{psiendpoint}
\frac{N^{c-1}\mathbf 1_{B_N}(\tau_2,\xi_2)\mathbf 1_{A_N}(\tau_1,\xi_1)}{N^{-\frac12}}
\lesssim
|\Psi(\tau,\xi,\tau_1,\xi_1)\mathbf 1_{B_N}(\tau_2,\xi_2)\mathbf 1_{A_N}(\tau_1,\xi_1)|.
\end{equation} 
Now, taking $f_N,g_N\!\in\!\mathcal S(\mathbb R^2)$ such that 
$\mathbf 1_{A_N}\le g_N$, $\mathbf 1_{B_N}\le f_N$, 
$\|g_N\|_{L^2}\lesssim\|\mathbf 1_{A_N}\|_{L^2}$, $\|f_N\|_{L^2}\lesssim\|\mathbf 1_{B_N}\|_{L^2}$
and combining \eqref{NE2eqv}, \eqref{psiendpoint} and \eqref{convbaixo}, follows the estimate
\begin{equation*}
N^{c-\frac12}
\lesssim 
1, 
\end{equation*}
which is false for $N$ sufficiently large whenever $c>1/2$. 
\hfill$\square$\\ 
%
%
%
\begin{ackn}
This paper is 
part of my Ph.D. thesis at the Federal University of Rio de Janeiro 
under the guidance of my advisor Didier Pilod. I want to take the opportunity 
to express my sincere gratitude to him. I also thank my colleagues at DMA/CEUNES 
in the Federal University of Espírito Santo for the support. 
The author was partially supported by CNPq-Brazil. 
\end{ackn} 
%
%
%
\indent 
 
%
%
\end{document}